\documentstyle[12pt,draft]{article}

 \newcounter{abceqn}

 \newcounter{abcfig}

\newcommand{\HH}{{\mathcal H}}
\newcommand{\LL}{{\mathcal L}}

\newcommand{\im}{\mathop{\rm Im}\nolimits}

\newcommand{\na}{\nabla}
\newcommand{\om}{\omega}
\newcommand{\Ga}{\Gamma}
\newcommand{\Om}{\Omega}

\newcommand{\Z}{Z^2/\{0\}}
\newcommand{\pa}{\partial}

\newcommand{\tz}{\tilde{z}}

\newcommand{\ta}{\tilde{a}}
\newcommand{\tN}{\tilde{N}}
\newcommand{\tw}{\tilde{w}}

\newcommand{\e}{\epsilon}

\newcommand{\k}{\kappa}
\newcommand{\ga}{\gamma}

\newcommand{\dl}{\delta}
\newcommand{\Dl}{\Delta}
\newcommand{\th}{\theta}
\newcommand{\ra}{\rightarrow}
\newcommand{\al}{\alpha}
\newcommand{\sg}{\sigma}
\newcommand{\Sg}{\Sigma}

\newcommand{\la}{\lambda}
\newcommand{\tla}{\tilde{\lambda}}

\newcommand{\nid}{\noindent}

\newcommand{\hk}{\hat{k}}

\renewcommand{\thesection}{\Roman{section}.}

\newcommand{\eqnsection}[1]{
	\section{#1}
	\setcounter{equation}{0}
	\renewcommand{\theequation}{\thesection\arabic{equation}}
	\setcounter{figure}{0}
	\renewcommand{\thefigure}{\arabic{figure}}
	\setcounter{remark}{0}
	\renewcommand{\theremark}{\thesection\arabic{remark}}
	\setcounter{theorem}{0}
	\renewcommand{\thetheorem}{\thesection\arabic{theorem}}
	\setcounter{lemma}{0}
	\renewcommand{\thelemma}{\thesection\arabic{lemma}}
}

\title{{\bf On 2D Euler Equations: \Large Part I.
On the Energy-Casimir Stabilities and The 
Spectra for Linearized 2D Euler Equations}}

\author{ \\ \\ \\ \\ 
Yanguang (Charles)\ \ Li  \thanks{This work is supported 
by the AMS Centennial Fellowship, and the Guggenheim Fellowship.}
\\  \\  \\ Department of Mathematics, 
 \\ \\ University of Missouri \\ \\ 
Columbia, MO 65211}

\date{\today}

\begin{document}
\bibliographystyle{unsrt}
\maketitle
\newpage
\begin{abstract}    
In this paper, we study a linearized two-dimensional Euler equation. 
This equation decouples into infinitely many invariant subsystems. Each 
invariant subsystem 
is shown to be a linear Hamiltonian system of infinite dimensions. 
Another important invariant besides the Hamiltonian for each invariant
subsystem is found, and is utilized to prove an ``unstable disk 
theorem'' through a simple Energy-Casimir argument \cite{HMRW85}. The eigenvalues of the 
linear Hamiltonian system are of four types: real pairs ($c,-c$), purely 
imaginary pairs ($id,-id$), quadruples ($\pm c\pm id$), and zero eigenvalues. 
The eigenvalues are computed through continued fractions.
The spectral equation for each invariant subsystem
is a Poincar\'{e}-type difference equation, i.e. it can be 
represented as the spectral equation of an infinite matrix operator, 
and the infinite matrix operator is a sum of a constant-coefficient 
infinite matrix operator and a compact infinite matrix operator. 
We have obtained a complete spectral theory.
\end{abstract}

\newtheorem{lemma}{Lemma}
\newtheorem{theorem}{Theorem}
\newtheorem{corollary}{Corollary}
\newtheorem{remark}{Remark}
\newtheorem{definition}{Definition}
\newtheorem{proposition}{Proposition}
\newtheorem{assumption}{Assumption}

\newpage
\setlength{\baselineskip}{28pt}

\eqnsection{Introduction}

In \cite{HK90}, Henshaw and Kreiss numerically studied the propagation 
of perturbations in the solution of the two dimensional incompressible 
Navier-Stokes equations with no body force, and at high Reynolds numbers. 
They numerically solved the equations starting with some smooth initial 
data for which the Fourier modes have random phases. This flow evolved 
over time, initially through a complicated state containing many shear 
layers and then into a state of large vortex structures (with shear 
layers) which persists for long time.
\begin{enumerate}
\item Changing the viscosity, the large vortex structures do not change 
much.
\item Adding high mode perturbation to the initial data, the large vortex 
structures do not change much.
\item Changing the Laplacian $\Dl$ to $\Dl^2$, the large vortex structures 
do not change much.
\item Changing the Laplacian $\Dl$ to an operator which is $\Dl$ for low 
modes, and is $\Dl^2$ for high modes, the large vortex structures do not 
change much.
\end{enumerate}
In \cite{MSMOM91}, Matthaeus et al. studied the same problem, and found 
similar results. Matthaeus et al. run the numerics for much longer time. 
These numerics indicate that for relatively long time (not infinite long 
time), the solution to the 2D N-S equation (without body force, i.e. 
decaying turbulence) has the large vortex structures. Such structures 
persist in the solution to the 2D Euler  equation by the claims 1, 3, 4 
above. (Cf: Under decay boundary conditions, the Kato's theorem states 
that for finite time the solution to 2D N-S equation converges to the 
solution to 2D Euler equation in norm as viscosity approaches zero, 
see \cite{Kat86}.) Moreover, the large vortex structures are stable 
with respect to the change in initial data.

In \cite{RS91}, Robert and Sommeria studied the organized structures in 
two-dimensional Euler fluid flows by a theory of equilibrium statistical 
mechanics. The theory takes into account all the known constants of motion 
for the two-dimensional Euler equations. The microscopic states are all 
the possible vorticity fields, while a macroscopic state is defined as a 
probability distribution of vorticity at each point of the domain, which 
descibes in a statistical sense the fine-scale vorticity fluctuations. 
The organized structure appears as a state of maximal entropy, with the 
constraints of all the constants of motion. The vorticity field obtained 
as the local average of this optimal macrostate is a {\em{steady solution}} 
of the Euler equations.

The above numerical results show that certain relatively long time large 
vortex structures for 2D Euler equation persist for 2D N-S equation at high 
Reynolds numbers. Such structures are stable with respect to the change 
in initial data. We believe such structures are the exhibition of certain 
unstable manifolds. The above theoretical results show that the probabiltity 
mean of such structures for 2D Euler flows are steady solutions to the 2D 
Euler equations. Thus, we believe that certain unstable manifolds, if 
there is any, of steady solutions to 2D Euler equations are responsible 
for such large vortex structures. Therefore, a dynamical system study 
on certain unstable manifolds of certain steady solutions to 2D Euler 
equations, and their persistence for 2D N-S equations at high Reynolds 
numbers is crucial for studying the large vortex structures. Indeed, we 
have built a {\em{chaos-molecules-model}} on 2D turbulence based upon the 
above motivations \cite{Li97b}. We believe that our chaos-molecules-model 
captures the qualititive frames of the hyperbolic structures, and the 
energy inverse cascade and the enstrophy cascade nature for 2D turbulence.

In this paper, we study the linearized 2D Euler equation at a fixed point. 
This study is the base for future analytical studies on the unstable 
manifolds for the 2D Euler equation. In \cite{WL97}, we have begun 
numerical studies on the unstable manifolds for the 2D Euler equation.

The current study is also important in the linear hydrodynamic stability 
theory. By utilizing Energy-Casimir method \cite{HMRW85}, we obtain an {\em{unstable 
disk theorem}} which is not in the category of the classical Rayleigh 
theorem \cite{Lin55}. 

Next we discuss the approaches used in this study. Through the 
Energy-Casimir method, nonlinear stabilities of various types of 
two-dimensional ideal fluid flows have been established \cite{Arn65} 
\cite{Arn69} \cite{Mar92} \cite{HMRW85}. Below we give a brief description on the 
Energy-Casimir method. Let $D$ be a region on the ($x,y$)-plane bounded 
by the curves $\Ga_i$ ($i=1,2$), an ideal fluid flow in $D$ is governed 
by the 2D Euler equation written in the stream-function form:
\begin{equation}
{\pa \over \pa t} \Dl \psi = [ \na \psi, \na \Dl \psi ]\ , \label{sfef}
\end{equation}
where 
\[
[ \na \psi, \na \Dl \psi ] = {\pa \psi \over \pa x}
{\pa \Dl \psi \over \pa y}
- {\pa \psi \over \pa y}{\pa \Dl \psi \over \pa x}\ ,
\]
with the boundary conditions,
\[
\psi|_{\Ga_i}=c_i(t)\ ,\ \ c_1 \equiv 0\ ,\ \ {d \over dt} \oint_{\Ga_i}
{\pa \psi \over \pa n} ds = 0\ .
\]
For every function $f(z)$, the functional 
\begin{equation}
F= \int\int_{D} f(\Dl \psi) \ dxdy \label{brs1}
\end{equation}
is a constant of motion (a Casimir) for (\ref{sfef}). The conditional 
extremum of the kinetic energy
\begin{equation}
E={1 \over 2}\int\int_{D}\na \psi \cdot \na \psi  \  dxdy \label{brs2}
\end{equation}
for fixed $F$ is given by the Lagrange's formula \cite{Arn65},
\begin{equation}
\dl H = \dl (E+\la F) =0\ ,\ \ \ \ \Rightarrow \ \ \psi_0 = 
\la f'(\Dl \psi_0)\ . 
\label{brs3}
\end{equation}
where $\la$ is the Lagrange multiplier. Thus, $\psi_0$ is the stream 
function of a stationary flow, which satisfies 
\begin{equation}
\psi_0 = \Phi(\Dl \psi_0)\ , \label{brs4}
\end{equation}
where $\Phi = \la f'$. The second variation is given by \cite{Arn65},
\begin{equation}
\dl^2 H = {1 \over 2} \int\int_{D} \bigg \{ \na \phi \cdot \na \phi + 
\Phi'(\Dl \psi_0) \ (\Dl \phi)^2 \bigg \} dxdy\ . \label{brs5}
\end{equation}
Let $\psi = \psi_0 + \varphi$ be a solution to the 2D Euler equation 
(\ref{sfef}), Arnold proved the estimates \cite{Arn69}: (a). when $c \leq 
\Phi'(\Dl \psi_0) \leq C$, $0 < c \leq C <\infty$,
\[
\int\int_{D} \bigg \{ \na \varphi(t) \cdot \na \varphi(t) + c 
[\Dl \varphi(t)]^2 \bigg \} \ dxdy
\leq \int\int_{D} \bigg \{ \na \varphi(0) \cdot \na \varphi(0) + C 
[\Dl \varphi(0)]^2 \bigg \} \ dxdy,
\]
for all $t \in (-\infty, +\infty)$, (b). when $c \leq 
-\Phi'(\Dl \psi_0) \leq C$, $0 < c < C <\infty$,
\[
\int\int_{D} \bigg \{ c [\Dl \varphi(t)]^2 - \na \varphi(t) \cdot 
\na \varphi(t)\bigg \} \ dxdy
\leq \int\int_{D} \bigg \{ C [\Dl \varphi(0)]^2 - \na \varphi(0) 
\cdot \na \varphi(0) \bigg \} \ dxdy,
\]
for all $t \in (-\infty, +\infty)$. Therefore, when the second variation 
(\ref{brs5}) is positive definite, or when 
\[
\int\int_{D} \bigg \{ \na \phi \cdot \na \phi + 
[\max \Phi'(\Dl \psi_0)] \ (\Dl \phi)^2 \bigg \} \ dxdy
\]
is negative definite, the stationary flow (\ref{brs4}) is nonlinearly stable 
(Liapunov stable). In this paper, we have found an invariant for the 
linearized 2D Euler equation, and use this invariant together with an 
Energy-Casimir type argument to study linear stability, and to prove an 
unstable 
disk theorem. The linearized 2D Euler equation is an infinite dimensional 
linear Hamiltonian system. 
For finite dimensional linear Hamiltonian systems, 
it is well-known that the eigenvalues are of 
four types: real pairs ($c, -c$), purely imaginary pairs ($id, -id$), 
quadruples ($\pm c \pm id$), and zero eigenvalues \cite{Poi99} \cite{Lia49} 
\cite{Arn80}. The same is true for the linearized 2D Euler equation. The 
eigenvalues are computed through continued fractions following the work of 
Meshalkin and Sinai \cite{MS61}.
The linearized 2D Euler equation can also be written in an 
infinite matrix form. The spectral equation of the infinite matrix operator 
defines a Poincar\'{e}-type difference equation \cite{Per10a} \cite{Per10b}. 
That is, the infinite matrix operator can be written as the sum of a 
constant-coefficient infinite matrix operator and a compact infinite matrix 
operator. In this paper, we follow a spectral theory developed by Duren 
\cite{Dur60} to study the spectra of the 
constant-coefficient 
infinite matrix operator through characteristic polynomials.
Then we apply the Weyl's essential spectra 
theorem to the perturbation of the constant-coefficient 
infinite matrix operator by the compact infinite matrix operator 
\cite{RS78}, to achieve a complete spectral theory. 

Finally, we discuss some preliminaries on the 2D Euler equation. Consider 
the two-dimensional incompressible Euler equation
written in vorticity form,
\begin{eqnarray}
& &{\pa \Om \over \pa t}= - u \ {\pa \Om \over \pa x} -v
\ {\pa \Om \over \pa y},
\nonumber \\
\label{Euler} \\
& &{\pa u \over \pa x} +{\pa v \over \pa y} = 0; \nonumber
\end{eqnarray}
under periodic boundary conditions in both $x$ and $y$ directions
with period $2\pi$, where $\Om$ is vorticity, $u$ and $v$ are 
respectively velocity components along $x$ and $y$ directions.
We also require that both 
$u$ and $v$ have means zero,
\[
\int_0^{2\pi}\int_0^{2\pi} u\ dxdy =\int_0^{2\pi}\int_0^{2\pi} v\ dxdy=0.
\]
Expand $\Om$ into Fourier series,
\[
\Om =\sum_{k\in Z^2/\{0\}} \om_k \ e^{ik\cdot X}\ ,
\]
where $\om_{-k}=\overline{\om_k}\ $, $k=(k_1,k_2)^T$, 
$X=(x,y)^T$. In this paper, we confuse $0$ with $(0,0)^T$, the context 
will always make it clear. By the relation between vorticity $\Om$ and stream 
function $\Psi$,
\[
\Om ={\pa v \over \pa x} -{\pa u \over \pa y} =\Dl \Psi \ ,
\]
where the stream function $\Psi$ is defined by,
\[
u=- {\pa \Psi \over \pa y} \ ,\ \ \ v={\pa \Psi \over \pa x} \ ;
\]
the system (\ref{Euler}) can be rewritten as the following kinetic system,
\begin{equation}
\dot{\om}_k = \sum_{k=p+q} A(p,q) \ \om_p \om_q \ ,
\label{Keuler}
\end{equation}
where $A(p,q)$ is given by,
\begin{eqnarray}
A(p,q)&=& {1\over 2}[|q|^{-2}-|p|^{-2}](p_1 q_2 -p_2 q_1) \nonumber \\
\label{Af} \\      
      &=& {1\over 2}[|q|^{-2}-|p|^{-2}]\left | \begin{array}{lr} 
p_1 & q_1 \\ p_2 & q_2 \\ \end{array} \right | \ , \nonumber
\end{eqnarray}
where $|q|^2 =q_1^2 +q_2^2$ for $q=(q_1,q_2)^T$, similarly for $p$.

\begin{remark}
Notice that direct calculation shows that the nonlinear term in (\ref{Keuler})
is $\sum_{k=p+q} \tilde{A}(p,q) \ \om_p \om_q$, where $\tilde{A}(p,q) = 
|p|^{-2}(q_1 p_2 - p_1 q_2)$; then, $\tilde{A}(q,p) = 
|q|^{-2}(p_1 q_2 - q_1 p_2)$. The $A(p,q)$ in (\ref{Keuler}) is the average 
of $\tilde{A}(p,q)$ and $\tilde{A}(q,p)$, $A(p,q)={1\over 2}[\tilde{A}(p,q) 
+ \tilde{A}(q,p)]$.
\end{remark}
For any two functionals $F_1$ and $F_2$ of $\{ \om_k \}$, 
define their Lie-Poisson bracket:
\begin{equation}
\{ F_1,F_2 \} = \sum_{k+p+q=0} \left | \begin{array}{lr} 
q_1 & p_1 \\ q_2 & p_2 \\ \end{array} \right | \ \om_k \
{\pa F_1 \over \pa \overline{\om_p}} \ {\pa F_2 \over \pa \overline{\om_q}}\ .
\label{Liebr}
\end{equation}
Then the 2D Euler equation (\ref{Keuler}) is a Hamiltonian system \cite{Arn66},
\begin{equation}
\dot{\om}_k = \{ \om_k, H\}, \label{hEft}
\end{equation}
where the Hamiltonian $H$ is the kinetic energy,
\begin{equation}
H= {1\over 2} \sum_{k \in Z^2/\{0\}} |k|^{-2} |\om_k |^2. 
\label{Ih}
\end{equation}
Following are Casimirs (i.e. invariants that Poisson commute with 
any functional) of the Hamiltonian system (\ref{hEft}):
\begin{equation}
J_n = \sum_{k_1 + \cdot \cdot \cdot +k_n =0} \om_{k_1} 
\cdot \cdot \cdot \om_{k_n}. \label{Ic}
\end{equation}
The following Proposition is concerned with the equilibrium manifolds 
of the 2D Euler equation (\ref{Keuler}).
\begin{proposition}
For any $k \in Z^2/\{0\}$, the infinite dimensional space
\[
E^1_k \equiv \bigg \{ \{ \om_{k'}\} \ \bigg | \ \om_{k'}=0,
\ \mbox{if}\ k' \neq rk,\ \forall r \in R \bigg \} \ ,
\]
and the finite dimensional space 
\[
E^2_k \equiv \bigg \{ \{ \om_{k'}\} \ \bigg | \ \om_{k'}=0,
\ \mbox{if}\ |k'| \neq |k| \bigg \} \ ,
\]
entirely consist of fixed points of the system (\ref{Keuler}).
\label{eman}
\end{proposition}

Proof: Let $k^0 \in \Z$, $\{ \om^0_k \} \in E^1_{k^0}$. For 
any $p,q \in \Z$, $\om^0_p \om^0_q \neq 0$ implies that $p = r_1 k^0,\ 
q=r_2 k^0$ for some $r_1, r_2 \in R$; then, $A(p,q)=0$; thus we 
always have $A(p,q)\ \om^0_p \om^0_q = 0$, and $\{ \om^0_k \}$ is a 
fixed point of (\ref{Keuler}). Let $k^0 \in \Z$, $\{ \om^0_k \} 
\in E^2_{k^0}$. 
For any $p,q \in \Z$, $\om^0_p \om^0_q \neq 0$ implies that 
$|p|=|k^0|,\ |q|=|k^0|$; then, $A(p,q)=0$; thus we 
always have $A(p,q)\ \om^0_p \om^0_q = 0$, and $\{ \om^0_k \}$ is a 
fixed point of (\ref{Keuler}). $\Box$

\begin{figure}[ht]
  \begin{center}
    \leavevmode
      \setlength{\unitlength}{2ex}
  \begin{picture}(36,27.8)(-18,-12)
    \thinlines
\multiput(-12,-11.5)(2,0){13}{\line(0,1){23}}
\multiput(-16,-10)(0,2){11}{\line(1,0){32}}
    \thicklines
\put(0,-14){\vector(0,1){28}}
\put(-18,0){\vector(1,0){36}}
\put(0,15){\makebox(0,0){$k_2$}}
\put(18.5,0){\makebox(0,0)[l]{$k_1$}}
\qbezier(-5.5,0)(-5.275,5.275)(0,5.5)
\qbezier(0,5.5)(5.275,5.275)(5.5,0)
\qbezier(5.5,0)(5.275,-5.275)(0,-5.5)
\qbezier(0,-5.5)(-5.275,-5.275)(-5.5,0)
    \thinlines
\put(4,4){\circle*{0.5}}
\put(0,0){\vector(1,1){3.7}}
\put(4,-4){\circle*{0.5}}
\put(-4,4){\circle*{0.5}}
\put(2,2){\circle*{0.5}}
\put(6,6){\circle*{0.5}}
\put(8,8){\circle*{0.5}}
\put(10,10){\circle*{0.5}}
\put(-2,-2){\circle*{0.5}}
\put(-4,-4){\circle*{0.5}}
\put(-6,-6){\circle*{0.5}}
\put(-8,-8){\circle*{0.5}}
\put(-10,-10){\circle*{0.5}}
\put(-10.6,-10.6){\line(1,1){21.5}}
\put(4.35,4.65){$k$}
\put(12.5,-12.5){\makebox(0,0)[t]{$|k'|=|k|$}}
\put(12.4,-12.4){\vector(-1,1){8.2}}
\put(14.5,-2.5){\makebox(0,0)[t]{$k'=rk$}}
\put(14.4,-2.4){\vector(-1,1){8.2}}
  \end{picture}
  \end{center}
\caption{An illustration on the locations of the modes ($k'=rk$) and 
($|k'|=|k|$) in the definitions of $E^1_k$ and $E^2_k$ 
(Proposition \ref{eman}).}
\label{eulman}
\end{figure}
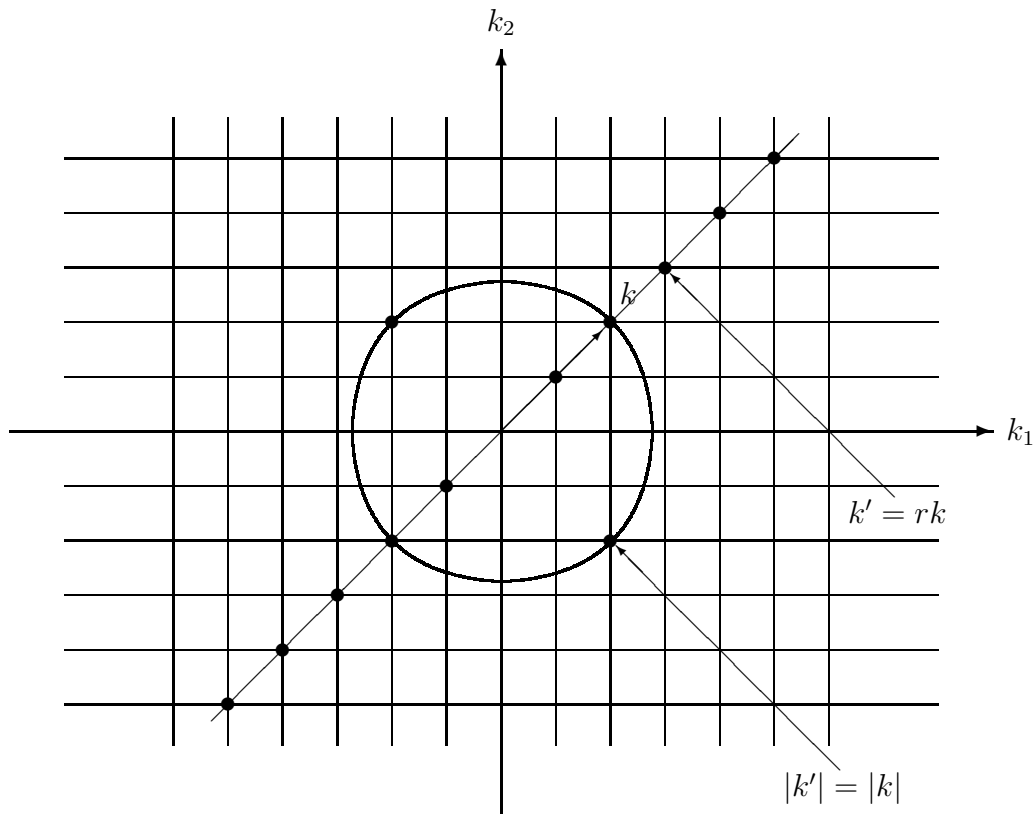

Fig.\ref{eulman} shows an example on the locations of the modes 
($k'=rk$) and ($|k'|=|k|$) in the definitions of $E_k^1$ and $E_k^2$
(Proposition \ref{eman}).

The paper is organized as follows: Section II is on the formulations of 
the problem. Section III is on the Liapunov stability. Section IV is on the
properties of the eigenvalues of the linearized 2D Euler equation as a 
linear Hamiltonian system. Section V is on the continued fraction study 
of the eigenvalues. Section VI is on the infinite-matrix study of 
the spectra of the linearized 2D Euler equation. Section VII is the 
conclusion.

\newpage
\eqnsection{The Formulations of the Problem}

Denote $\{ \om_k \}_{k\in \Z}$ by $\om$. Consider the simple fixed point 
$\om^*$:
\begin{equation}
\om^*_p = \Ga,\ \ \ \om^*_k = 0 ,\ \mbox{if} \ k \neq p \ \mbox{or}\ -p,
\label{fixpt}
\end{equation}
of the 2D Euler equation (\ref{Keuler}), which belongs to the two-dimensional 
intersection space $E^1_p \cap E^2_p$ (Proposition \ref{eman}), where 
$\Ga$ is an arbitrary complex constant. The 
{\em{linearized two-dimensional Euler equation}} at $\om^*$ is given by,
\begin{equation}
\dot{\om}_k = A(p,k-p)\ \Ga \ \om_{k-p} + A(-p,k+p)\ \bar{\Ga}\ \om_{k+p}\ .
\label{LE}
\end{equation}
This is the linearized two-dimensional Euler equation that we are going 
to study in this paper.

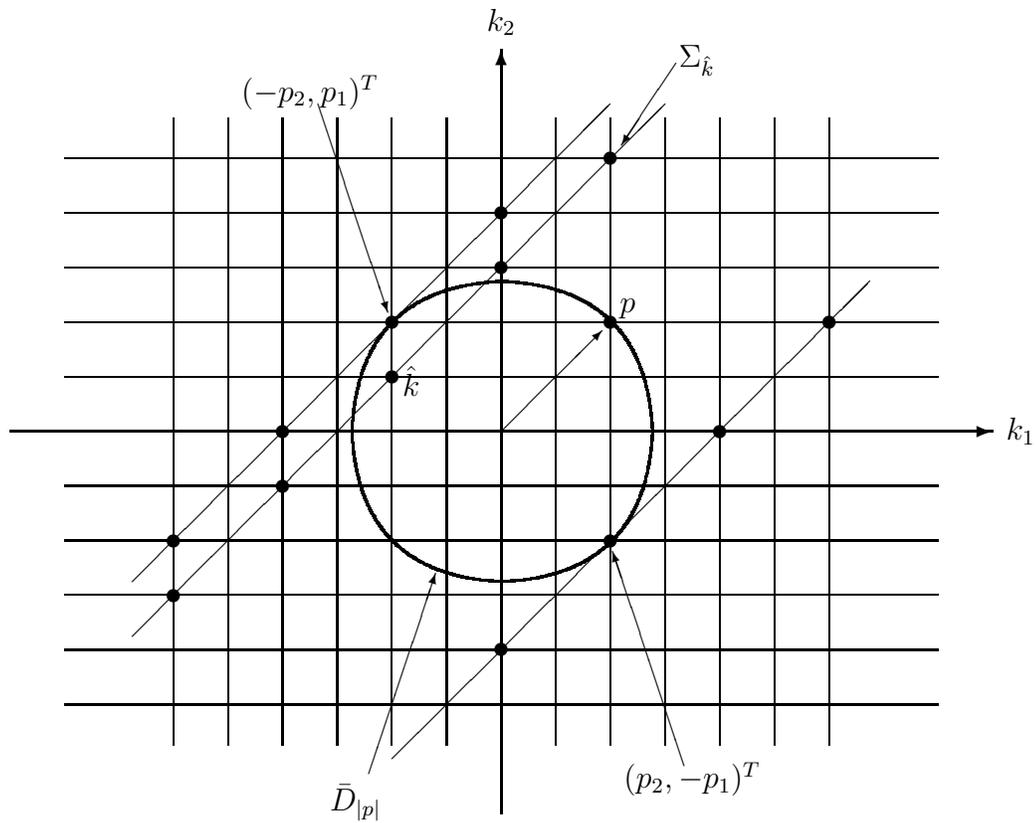
\begin{figure}[ht]
  \begin{center}
    \leavevmode
      \setlength{\unitlength}{2ex}
  \begin{picture}(36,27.8)(-18,-12)
    \thinlines
\multiput(-12,-11.5)(2,0){13}{\line(0,1){23}}
\multiput(-16,-10)(0,2){11}{\line(1,0){32}}
    \thicklines
\put(0,-14){\vector(0,1){28}}
\put(-18,0){\vector(1,0){36}}
\put(0,15){\makebox(0,0){$k_2$}}
\put(18.5,0){\makebox(0,0)[l]{$k_1$}}
%
%
%
%
%
\qbezier(-5.5,0)(-5.275,5.275)(0,5.5)
\qbezier(0,5.5)(5.275,5.275)(5.5,0)
\qbezier(5.5,0)(5.275,-5.275)(0,-5.5)
\qbezier(0,-5.5)(-5.275,-5.275)(-5.5,0)
    \thinlines
\put(4,4){\circle*{0.5}}
\put(0,0){\vector(1,1){3.7}}
\put(4.35,4.35){$p$}
\put(4,-4){\circle*{0.5}}
\put(8,0){\circle*{0.5}}
\put(-8,0){\circle*{0.5}}
\put(-8,-2){\circle*{0.5}}
\put(-12,-4){\circle*{0.5}}
\put(-12,-6){\circle*{0.5}}
\put(-4,2){\circle*{0.5}}
\put(-4,4){\circle*{0.5}}
\put(0,6){\circle*{0.5}}
\put(0,8){\circle*{0.5}}
\put(4,10){\circle*{0.5}}
\put(12,4){\circle*{0.5}}
\put(0,-8){\circle*{0.5}}
\put(-4,-12){\line(1,1){17.5}}
\put(-13.5,-7.5){\line(1,1){19.5}}
\put(-13.5,-5.5){\line(1,1){17.5}}
\put(-3.6,1.3){$\hat{k}$}
\put(-7,12.1){\makebox(0,0)[b]{$(-p_2, p_1)^T$}}
\put(-6.7,12){\vector(1,-3){2.55}}
\put(6.5,13.6){\makebox(0,0)[l]{$\Sg_{\hat{k}}$}}
\put(6.4,13.5){\vector(-2,-3){2.0}}
\put(7,-12.1){\makebox(0,0)[t]{$(p_2, -p_1)^T$}}
\put(6.7,-12.25){\vector(-1,3){2.62}}
\put(-4.4,-13.6){\makebox(0,0)[r]{$\bar{D}_{|p|}$}}
\put(-4.85,-12.55){\vector(1,3){2.45}}
  \end{picture}
  \end{center}
\caption{An illustration of the classes $\Sg_{\hk}$ and the disk 
$\bar{D}_{|p|}$.}
\label{class}
\end{figure}

\begin{definition}[Classes]
For any $\hk \in \Z$, we define the class $\Sg_{\hk}$ to be the subset of 
$\Z$:
\[
\Sg_{\hk} = \bigg \{ \hk + n p \in \Z \ \bigg | \ n \in Z, \ \ p \ \mbox{is 
specified in (\ref{fixpt})} \bigg \}.
\]
\label{classify}
\end{definition}
See Fig.\ref{class} for an illustration of the classes. 
According to the classification 
defined in Definition \ref{classify}, the linearized two-dimensional Euler 
equation (\ref{LE}) decouples into infinite many invariant subsystems:
\begin{eqnarray}
\dot{\omega}_{\hat{k} + np} &=& A(p, \hat{k} + (n-1) p) 
     \ \Gamma \ \omega_{\hat{k} + (n-1) p} \nonumber \\  \label{CLE}\\
& & + \ A(-p, \hat{k} + (n+1)p)\ 
     \bar{\Gamma} \ \omega_{\hat{k} +(n+1)p}\ . \nonumber
\end{eqnarray}
Each invariant subsystem can be rewritten as a linear Hamiltonian system 
as shown below.
\begin{definition}[The Quadratic Hamiltonian] 
The quadratic Hamiltonian $\HH_{\hat{k}}$ is defined as:
\begin{eqnarray}
\HH_{\hat{k}} &=& -2 \ \im \bigg \{ \sum_{n \in Z} \rho_n \ \Gamma 
\ A(p, \hat{k} + (n-1)p)\ \omega_{\hat{k} + (n-1)p} \ 
\bar{\omega}_{\hat{k} +np} \bigg \} \nonumber \\
\label{CHAM} \\             
&=& - \left| 
      \begin{array}{cc}
          p_1 & \hat{k}_1\\
          p_2 & \hat{k}_2
      \end{array}
      \right| \ \im
      \bigg \{ \sum_{n \in Z} \Gamma \ \rho_n \ \rho_{n-1}
        \ \omega_{\hat{k} + (n-1)p}
        \ \bar{\omega}_{\hat{k} + np} 
      \bigg \}, \nonumber
\end{eqnarray}
where $\rho_n = [ |\hat{k} + np |^{-2}-|p|^{-2}]$,
``\ $\im$\ '' denotes `` imaginary part ''. 
\end{definition}
Then the invariant subsystem (\ref{CLE}) can be rewritten as a linear 
Hamiltonian system,
\begin{equation}
i \ \dot{\omega}_{\hat{k} +n p} = \rho^{-1}_n \ \frac{\pa \HH_{\hat{k}}}
         {\partial \bar{\omega}_{\hat{k} + np}}\ .
\label{CHAF}
\end{equation}
Let
\[
  \omega_{\hat{k} +np} = \alpha_n + i \beta_n \, , 
    \quad n \in Z \, ,
\]
i.e. $\alpha_n$ and $\beta_n$ are the real and imaginary parts of
$\omega_{\hat{k} +np}$. Then the linear Hamiltonian system (\ref{CHAF}) 
can be rewritten in the form,
\begin{equation}
   \left\{
    \begin{array}{rcl}
      \dot{\alpha}_n &=& \frac{1}{2} \ \rho^{-1}_n \
          \frac{\partial \HH_{\hat{k}}}
          {\partial \beta_n} \, , \\[2ex]
& & \\ 
      \dot{\beta}_n &=& - \frac{1}{2}\  \rho^{-1}_n 
          \ \frac{\partial \HH_{\hat{k}}}{\partial \alpha_n} \, ,
    \end{array}
   \right.
\label{LHR}
\end{equation}
where
\begin{eqnarray*}
  \HH_{\hat{k}} &=& \left| 
    \begin{array}{cc}
      \hat{k}_1 & p_1\\
      \hat{k}_2 & p_2
    \end{array}\right| \
    \sum_{n \in Z} \rho_n \ \rho_{n-1}\  
    \bigg [ \Gamma_r \ ( \alpha_n \beta_{n-1} -
      \alpha_{n-1} \beta_n) \nonumber \\ \nonumber \\
& & + \ \Gamma_i \ (\alpha_{n-1} \alpha_n + \beta_{n-1} \beta_n) 
\bigg ] \, , \nonumber
\end{eqnarray*}
where  $\Gamma = \Gamma_r + i \ \Gamma_i$, i.e. $ \Gamma_r$ and
$\Gamma_i$ are the real and imaginary parts of $\Gamma$. Explicitly,
\begin{equation}
  \left\{
    \begin{array}{rcl}
      \dot{\alpha}_n &=& \frac{1}{2}
        \left|
         \begin{array}{cc}
           \hat{k}_1 & p_1\\
           \hat{k}_2 & p_2
         \end{array}
        \right|
        \ \bigg [ \rho_{n+1} \ \Gamma_r \ \alpha_{n+1} -
              \rho_{n-1} \ \Gamma_r \ \alpha_{n-1} \\  \\
            & & +
              \rho_{n+1} \ \Gamma_i \ \beta_{n+1} +
              \rho_{n-1} \ \Gamma_i \ \beta_{n-1} 
        \bigg ] \, , \\  \\
      \dot{\beta}_n &=& - \frac{1}{2}
        \left|
          \begin{array}{cc}
            \hat{k}_1 & p_1 \\
            \hat{k}_2 & p_2
          \end{array}
        \right|
        \ \bigg [ \rho_{n-1} \ \Gamma_r \ \beta_{n-1} -
               \rho_{n+1} \ \Gamma_r \ \beta_{n+1} \\  \\
           & & +
               \rho_{n+1} \ \Gamma_i \ \alpha_{n+1} +
               \rho_{n-1} \ \Gamma_i \ \alpha_{n-1}
        \bigg ] \, .
    \end{array}
    \right.
\label{rlh}
\end{equation}
If we rescale the variables as follows,
\[
  \tau = \frac{1}{2}
     \left|
       \begin{array}{cc}
         \hat{k}_1 &  p_1 \\
         \hat{k}_2 &  p_2
       \end{array}
     \right| \ t \ ,
\ \ \tilde{\alpha}_n = \rho_n  \ \alpha_n  \, , \ \ 
     \tilde{\beta}_n = \rho_n \ \beta_n \, ,
\]
The linear Hamiltonian system (\ref{rlh}) can be rewritten in the simpler 
form,
\[
  \left\{
    \begin{array}{rcl}
      d \tilde{\alpha}_n / d \tau &=& \rho_n 
        \ \bigg [ 
          \Gamma_r \ (\tilde{\alpha}_{n+1} - \tilde{\alpha}_{n-1})\ + \ 
          \Gamma_i \ ( \tilde{\beta}_{n+1} + \tilde{\beta}_{n-1} )
        \bigg ] \, , \\  \\ 
     d \tilde{\beta}_n / d \tau &=& - \rho_n
        \ \bigg [ 
          \Gamma_r \ (\tilde{\beta}_{n-1} - \tilde{\beta}_{n+1}) \ + \ 
          \Gamma_i \ (\tilde{\alpha}_{n+1} + \tilde{\alpha}_{n-1})
        \bigg ]\ .
    \end{array}
    \right. 
\]
Next we discuss the fixed point $\om^*$ from a variational-principle 
point of view. Consider the kinetic energy and the enstrophy of the 
2D Euler equation (\ref{Keuler}),
\[
E = \frac{1}{2} \sum\limits_{k \in Z^2 / \left\{ 0 \right\}}
        \left| k \right|^{-2} 
        \left| \omega_k \right|^2 \ , \ \ \ 
J = \sum\limits_{k \in Z^2 / \left\{ 0 \right\}}
        \left| \omega_k \right|^2 \ .
\]
If we extremize the kinetic energy $E$ for fixed enstrophy $J=c$ (a constant), 
we have the critical states by the Lagrange formula,
\[
\left\{
    \begin{array}{rcl}
      \frac{\partial L}{\partial \lambda} &=& J-c = 0 \, , \\ \\
      \frac{\partial L}{\partial \omega_k} &=&
          (\left| k \right|^{-2} + 
          \lambda)\ \bar{\omega}_k = 0,  \quad 
       \forall k \in Z^2 / \left\{ 0 \right\}\, ,
    \end{array}
  \right. 
\]
where $L = 2 E + \lambda (J-c)$. Thus, the critical states satisfy the 
relations,
\begin{equation}
\left \{  \begin{array}{l}
    \lambda = -\left | q \right |^{-2} \, , \quad \mbox{for some}\ 
q \in \Z, \\ \\ 
       \omega_k =0 \, , \quad 
       \hbox{if } \left | k \right | \ne \left | q \right | \, , \\ \\
    \sum\limits_{\left | k \right | = \left | q \right |} \left | \omega_k
       \right |^2 = c \, . \\
\end{array} \right.
\nonumber
\end{equation}
Thus, we have the {\em{critical manifold}} for extremizing the kinetic 
energy for fixed enstrophy,
\[
M^c_{\left| q \right|} =
  \bigg \{ \omega\ \bigg |\ \omega_k = 0,\ \hbox{ if }\ \left| k \right| 
\ne \left| q \right| , \ \ \sum_{\left| k \right| = \left| q \right|}
    \left| \omega_k \right|^2 = c \bigg \},
\]
which is a submanifold of the equilibrium manifold $E^2_q$ (Proposition 
\ref{eman}). Denote by $I$ the linear combination,
\begin{equation}
I =  2E - \left| p \right|^{-2} J.
\label{fuI}
\end{equation}
Then, we have the following variational principle for 
the fixed point $\om^*$ (\ref{fixpt}). 

{\bf Variational Principle:} The fixed point
$\omega^*$ is a conditionally critical state of the kinetic energy $E$ for
fixed enstrophy $J = | \Gamma |^2$, and is an
absolute critical state of $I$. 

\eqnsection{Liapunov Stability}

\begin{definition}[An Important Functional]
For each invariant subsystem (\ref{CLE}), we define the functional 
$I_{\hat{k}}$
which is the restriction of the functional $I$ (\ref{fuI}) to the class 
$\Sg_{\hat{k}}$,
\begin{eqnarray}
    I_{\hat{k}} &=& I_{(\hbox{\textup{\small restricted to }} 
                       \Sg_{\hat{k}})} \nonumber \\
\label{invIk}   \\            
&=& \sum_{n \in Z} \{ | \hat{k} + np |^{-2}
                    - | p |^{-2} \} \left|
                    \omega_{\hat{k} + np} \right|^2 \, . \nonumber
\end{eqnarray}
\end{definition}

\begin{lemma}
$I_{\hat{k}}$ is a constant of motion for the system (\ref{CHAF}).
\label{const}
\end{lemma}

Proof: Differentiating $I_{\hat{k}}$, we have
\begin{eqnarray*}
    \dot{I}_{\hat{k}} &=& \sum_{n \in Z} \rho_n
        \left[ \dot{\omega}_{\hat{k} + np} \bar{\omega}_{\hat{k} +np}
             + \omega_{\hat{k} +np} \dot{\bar{\omega}}_{\hat{k} +np}
        \right] \\ \\
     &=& -i\ \sum_{n \in Z} 
         \left[ \bar{\omega}_{\hat{k} +np}\ 
         \frac{\partial \HH_{\hat{k}}}
              {\partial \bar{\omega}_{\hat{k} +np}} \ - \ 
              \omega_{\hat{k} + np} 
         \ \frac{\partial \HH_{\hat{k}}}
              {\partial \omega_{\hat{k} + np}}\ 
          \right]\\ \\
      &=&  \frac{1}{2} 
           \left|
               \begin{array}{cc}
                 p_1 & \hat{k}_1\\
                 p_2 & \hat{k}_2
               \end{array}
           \right| \ 
           \sum_{n \in Z} 
           \bigg[ \bar{\omega}_{\hat{k} + np} 
              \bigg ( \Gamma \rho_n \rho_{n-1} \omega_{\hat{k} + (n-1) p} \\ \\
           & & - \bar{\Gamma} \rho_{n+1} \rho_n 
                  \omega_{\hat{k} + (n+1)p}
              \bigg ) 
- \omega_{\hat{k} +np} 
             \bigg ( \Gamma\rho_{n+1} \rho_n 
\bar{\omega}_{\hat{k} +(n+1)p} \\ \\ 
& & - \bar{\Gamma} \rho_n \rho_{n-1} 
                \bar{\omega}_{\hat{k} +(n-1)p} 
             \bigg ) 
         \bigg ] = 0 \, .
\end{eqnarray*}
This completes the proof of the lemma. $\Box$ 

Next we define the concept of disk in $\Z$ which is needed in the 
unstable disk theorem to be proved below.
\begin{definition}[The Disk]
The disk of radius $\left| p \right|$ in $Z^2 / \left\{ 0
\right\}$, denoted by $D_{\left| p \right|}$, is defined as
\[
  D_{\left| p \right|} = \bigg \{ k \in Z^2 / \left\{ 0 \right\} \ \bigg| 
      \ \left| k \right| < \left| p \right| \bigg \} \, .
\]
The closure of $D_{\left| p \right|}$, denoted by
$\bar{D}_{\left| p \right|}$, is defined as
\[ 
 \bar{D}_{\left| p \right|} = \bigg \{ k \in Z^2/ \left\{ 0 \right\} \ \bigg| 
     \ \left| k \right| \leq \left| p \right| \bigg \} \, .
\]
\end{definition}
See Fig.\ref{class} for an illustration.
Next we prove the unstable disk theorem using a simple Energy-Casimir type 
argument \cite{Arn66} \cite{HMRW85}.
\begin{theorem}[Unstable Disk Theorem]
If $\ \Sigma_{\hat{k}} \cap \bar{D}_{\left| p \right|} = \emptyset,\ $
then the invariant subsystem (\ref{CLE}) is Liapunov stable
for all $t \in R$, in fact, 
\[
  \sum_{n \in Z} \left|
  \omega_{\hat{k}+np}(t) \right|^2 \leq \sigma \ \sum_{n \in Z} \left|
  \omega_{\hat{k}+np}(0) \right|^2 \, , \quad \quad \forall t \in R \, , 
\]
where 
\[
  \sigma = \left[ \max_{n \in Z} 
              \left\{ - \rho_n \right\}
           \right] \, 
           \left[ \min_{n \in Z} 
               \left\{ -\rho_n \right\}
           \right]^{-1} \, , \quad 0< \sigma < \infty \, .
\]
\label{UDT}
\end{theorem}

Proof:
By lemma \ref{const}, $I_{\hat{k}}$ is a constant
of motion for the invariant subsystem (\ref{CLE}); then
\begin{equation}
    \sum_{n \in Z} \rho_n \left| \omega_{\hat{k}+np}(t) \right|^2 
      =  \sum_{n \in Z} \rho_n \left| \omega_{\hat{k} +np}(0)
         \right|^2 \, , 
         \quad \forall t \in R \, .
\label{plia0}
\end{equation}
If $\ \Sg_{\hat{k}} \cap \bar{D}_{\left| p \right|} = \emptyset$, then
\[
  | \hat{k} + np | > |p| \, , \quad
      \forall n \in Z \, . 
\]
Thus, there exists a constant $\delta >0$, such that
\begin{equation}
  \delta< - \rho_n < 2 \, .
\label{plia1}
\end{equation}
By (\ref{plia0}),
\begin{eqnarray*}
  \min_{n \in Z} \left\{ - \rho_n \right\} \sum_{n \in Z} 
     \left| \omega_{\hat{k} +np}(t) \right|^2
     &\leq& - \sum_{n \in Z} \rho_n \left| \omega_{\hat{k} +
              np}(t) \right|^2 \\
\\
     &=& - \sum_{n \in Z} \rho_n \left| \omega_{\hat{k}+np}(0)
            \right|^2  \\ \\
     &\leq& \max_{n \in Z} \left\{ - \rho_n \right\}
            \sum_{n \in Z} \left| \omega_{\hat{k}+np}(0) \right|^2 \, ,
\end{eqnarray*}
that is,
\[
\sum_{n \in Z} \left| \omega_{\hat{k}+np}(t) \right|^2 \leq
  \sigma \sum_{n \in Z} \left| \omega_{\hat{k}+np}(0) \right|^2 \ ,
\]
where 
\[
\sigma = \left[ \max_{n \in Z} \left\{ -\rho_n
      \right\} \right] \left[ \min_{n \in Z} \left\{ -\rho_n
      \right\} \right]^{-1}\ .  
\]
By relation (\ref{plia1}), 
\begin{displaymath}
      \frac{1}{2} \delta < \sigma < 2 \delta^{-1} \, .
\end{displaymath}
This completes the proof of the theorem. $\Box$
\begin{remark}
If $\ \hat{k} \parallel p$, i.e.\ $\exists$ real scalar $\alpha$ such
that $\ \hat{k} = \alpha \ p\ $, then the invariant subsystem (\ref{CLE})  
reduces to
\begin{displaymath}
    \dot{\omega}_{\hat{k}+np} = 0 \, , \quad \quad \forall n \in Z \, ;
\end{displaymath}
thus, it is obviously Liapunov stable for all $t \in R$, in fact,
this is a linearization inside the equilibrium space $E^1_p$ of the
$2D$ Euler equation (cf: Definition \ref{eman}).
\end{remark}
If $\ | \hat{k}| = |p|$, then the invariant subsystem
(\ref{CLE}) decomposes into two decoupled systems,
\begin{eqnarray}
\dot{\omega}_{\hat{k}+np} &=& A(p, \hat{k}+(n-1)p) \ \Gamma
  \omega_{\hat{k}+(n-1)p} \nonumber \\
\label{CLE1} \\
& & + A(-p, \hat{k}+(n+1)p) \ \bar{\Gamma}
  \omega_{\hat{k}+(n+1)p} \, , \quad (n \geq 1), \nonumber \\
\nonumber \\
  \dot{\omega}_{\hat{k} +np} &=& A(p, \hat{k} +(n-1)p) \
       \Gamma \omega_{\hat{k}+(n-1)p} \nonumber \\
\label{CLE2} \\
& &+ A(-p, \hat{k} +(n+1)p) \ 
       \bar{\Gamma} \omega_{\hat{k} +(n+1)p} \, , \quad (n \leq -1),
\nonumber
\end{eqnarray}
where $A(p, \hat{k})=0$.  The equation for $\omega_{\hat{k}}$ is
\begin{equation}
  \dot{\omega}_{\hat{k}} = A(p, \hat{k} -p) \ \Gamma \omega_{\hat{k}-p} +
                 A(-p, \hat{k} +p) \ \bar{\Gamma}
                 \omega_{\hat{k}+p} \, .
\label{CLE3}
\end{equation}
Each of (\ref{CLE1}) and (\ref{CLE2}) is a Hamiltonian system with
the Hamiltonian, \\
\begin{eqnarray*}
  \HH^+_{\hat{k}} &=& -      
  \left|
    \begin{array}{cc}
         p_1 & \hat{k}_1\\
         p_2 & \hat{k}_2  
    \end{array}
  \right| \ 
\im
   \left\{ \sum^{\infty}_{n=1} \Gamma \rho_n \rho_{n-1}
         \omega_{\hat{k}+(n-1)p} \bar{\omega}_{\hat{k}+np}
   \right\} \, , \\
\\ 
  \HH^-_{\hat{k}} &=& -      
  \left|
    \begin{array}{cc}
         p_1 & \hat{k}_1\\
         p_2 & \hat{k}_2  
    \end{array}
  \right| \ 
\im
   \left\{ \sum^{- \infty}_{n=-1} \Gamma \rho_n \rho_{n-1}
         \omega_{\hat{k}+(n-1)p} \bar{\omega}_{\hat{k}+np}
   \right\} \, ;
\end{eqnarray*}
which has the same representation as (\ref{CHAF}).
Denote by $I^+_{\hat{k}}$ and  $I^-_{\hat{k}}$ the restrictions of
$I_{\hat{k}}$ to the systems (\ref{CLE1}) and (\ref{CLE2}):
\begin{equation}
       I^+_{\hat{k}}= \sum^{\infty}_{n=1} \rho_n \left|
           \omega_{\hat{k}+np} \right|^2 \, ,
\label{rI1}
\end{equation}
\begin{equation}
       I^-_{\hat{k}}= \sum^{- \infty}_{n=-1} \rho_n \left|
           \omega_{\hat{k}+np} \right|^2 \, .
\label{rI2}
\end{equation}

\begin{lemma}
      If $\ | \hat{k}| = |p|$, then
      $I^+_{\hat{k}}$ and   $I^-_{\hat{k}}$ are respectively
      constants of motion for the systems (\ref{CLE1}) and (\ref{CLE2}).
\end{lemma}

Proof:
      The same with that for Lemma~\ref{const} . $\Box$
 
\begin{theorem}[Half Class Stability Theorem]
If $|\hat{k}| = |p|$, $\hat{k} +p \notin \bar{D}_{\left| p \right|}$, 
then the linear
Hamiltonian system (\ref{CLE1}) is Liapunov stable for all
$t \in R$, in fact,
\[
\sum^{\infty}_{n=1} \left| \omega_{\hat{k}+np}(t)
           \right|^2 \leq \sigma \sum^{\infty}_{n=1} \left|
           \omega_{\hat{k}+np}(0) \right|^2 \, , 
         \quad 
         \forall t \in R \, ,
\]
where 
\[
\sigma = \left[ \max_{n \geq 1} \left\{
           -\rho_n \right\} \right] \left[ \min_{n \geq  1} 
           \left\{ -\rho_n \right\} \right]^{-1} \, ,\quad 
             0<\sigma < \infty \ .
\]
If $ |\hat{k}| = | p |$, $\hat{k}
-p \notin \bar{D}_{\left| p \right|}$, then the linear
Hamiltonian system (\ref{CLE2}) is Liapunov stable for all
$t \in R$, in fact,
\[
\sum^{-\infty}_{n=-1} \left| \omega_{\hat{k}+np}(t)
           \right|^2 \leq \sigma \sum^{-\infty}_{n=-1} \left|
           \omega_{\hat{k}+np}(0) \right|^2 \, , 
         \quad 
         \forall t \in R \, ,
\]
where 
\[
\sigma = \left[ \max_{n \leq -1} \left\{
           -\rho_n \right\} \right] \left[ \min_{n \leq  -1} 
           \left\{ -\rho_n \right\} \right]^{-1} \, ,\quad 
             0<\sigma < \infty \ .
\]
\label{HLST}
\end{theorem}

    Proof:
      The same argument on $I^+_{\hat{k}}$ and  $I^-_{\hat{k}}$
      as that on  $I_{\hat{k}}$ in the proof of Theorem~\ref{UDT}.
    $\Box$

    \begin{remark}
      If $\hat{k} = (p_2, -p_1)^T$ or  $\hat{k} = (-p_2, p_1)^T$,
      then $\hat{k}+p \notin \bar{D}_{\left| p \right|}$ and
      $\hat{k}-p \notin \bar{D}_{\left| p \right|}$.  Therefore,
      by Theorem~\ref{HLST}, both systems (\ref{CLE1}) and
      (\ref{CLE2}) are Liapunov stable for all $t \in R$.  Points
      in $\sum_{\hat{k}}$ are on the tangent lines to the circle
        of radius $\left| p \right|$ at $\hat{k}$ (cf: Fig.\ref{class}).
    \end{remark}

\newpage
\eqnsection{Properties of the Point Spectrum for the Linearized
Two-Dimensional Euler Equation as a Linear Hamiltonian System}

In this section, we study the properties of the eigenvalues for the 
linear Hamiltonian system (\ref{LHR}). The right hand side of 
(\ref{LHR}) defines a
linear operator denoted by $\LL$, i.e.\
\begin{eqnarray}
  \LL \left( 
    \begin{array}{cc}
      \alpha \\
      \beta
    \end{array}
 \right) = 
     \left(
    \begin{array}{cc}
      \xi \\
       \eta 
    \end{array}
\right) \, ,
\label{spl}
\end{eqnarray}
where
\begin{eqnarray*}
  \begin{array}{rcl@{\qquad}rcl}
\alpha &=& (\cdots \alpha_{-1}\ \alpha_0 \ \alpha_1 \cdots )^T \, , 
        & \beta &=&  (\cdots \beta_{-1} \ \beta_0 \ \beta_1 \cdots )^T \, ; 
\\ \\ 
\xi &=& (\cdots \xi_{-1} \ \xi_0 \ \xi_1 \cdots )^T \, ,
         &  \eta &=&  (\cdots \eta_{-1} \ \eta_0 \ \eta_1 \cdots )^T \, ; 
\\ \\
\xi_n &=&  \frac{1}{2} \ \rho^{-1}_n \ 
           \frac{\partial \HH_{\hat{k}}}{\partial \beta_n} \, , 
      & \eta_n &=& - \frac{1}{2}\ \rho^{-1}_n \ 
            \frac{\partial \HH_{\hat{k}}}{\partial \al_n} \, , \quad
            n \in Z \, .
  \end{array}
\end{eqnarray*}
$\LL$ has the infinite matrix representation:  
\begin{eqnarray*}
\k \left(
  \begin{array}{c@{\!\!}ccc@{\!\!}c@{\;|\;}c@{\!\!}ccc@{\!\!}c}
           n \to - \infty & & & & & n \to - \infty &&&& \\[-0.2ex]
           \ddots \ \  & \ddots \ \ & \ddots \ \ &&& \ddots \ \ & 
          \ddots \ \ &\ddots \ \ && \\  
          &&&& \bigcirc &&&&& \bigcirc \\
&-u_{n-1}& \; 0 \; & u_{n+1} & &&
           v_{n-1}& \; 0 \; & v_{n+1} &  \\[-0.2ex]
\bigcirc  &&&&& \bigcirc  &&&& \\
&& \ \ \ddots & \ \ \ddots & \ \
\ddots & & &\ \ \ddots &\ \ \ddots &\ \ \ddots  \\ 
& & & n \to + \infty && & & &  n 
\to + \infty& \\
\hline
           n \to - \infty & & & & & n \to - \infty & & & & \\[-0.2ex]
           \ddots \ \ & \ddots \ \ &\ddots  \ \ & && \ddots \ \  &
\ddots \ \  &\ddots \ \ &&  \\
&&&&\bigcirc &&&&& \bigcirc \\
           & -v_{n-1}& \; 0 \; & -v_{n+1} && & 
           -u_{n-1}& \; 0 \; & u_{n+1}& \\[-0.2ex]
           \bigcirc  &&&&& \bigcirc &&&& \\
&&
\ \ \ddots & \ \ \ddots &\ \ \ddots  
& & &\ \ \ddots & \ \ \ddots &\ \ \ddots \\ 
& & & n \to + \infty & & & &&  n 
\to + \infty \\
\end{array}
\right)\!\!\!\!
\end{eqnarray*}
where $u_n = \Ga_r  \rho_n$, $v_n = \Ga_i \rho_n$, 
$\k = \frac{1}{2}\left| \begin{array}{cc}
                 \hat{k}_1 & p_1 \\
                 \hat{k}_2 & p_2
\end{array}\right| \ $.
We define the enstrophy norm which is the $\ell_2$ norm:
\begin{equation}
\bigg  \| ( \alpha , \beta)^T \bigg \|^2 = \sum_{n \in Z} \left( 
\alpha_n^2 + \beta_n^2 \right) \, .
\label{norm}
\end{equation}
\begin{lemma}
  The linear operator $\LL$ maps $\ell_2 \times \ell_2$ into
  $\ell_2 \times \ell_2 $:
\begin{displaymath}
    \LL \ :\ \ell_2 \times \ell_2 \mapsto \ell_2 \times \ell_2 \, .
  \end{displaymath}
\end{lemma}

Proof:
  Notice that
  \begin{displaymath}
    \rho_n \longrightarrow \left| p \right|^{-2} \, , 
       \quad  \hbox{as } \left| n \right| \longrightarrow \infty \, .
  \end{displaymath}
Let
\begin{displaymath}
  \rho_* = \max_{n \in Z} \left| \rho_n \right| 
\end{displaymath}
be the maximum of $\left| \rho_n \right|$; then $\rho_* < \infty$, and from
(\ref{spl}), we have
\begin{eqnarray*}
  \left| \xi_n \right| &\leq&  c \ 
        \bigg \{ \left| \alpha_{n+1} \right| +  \left|
                       \alpha_{n-1} \right| +
                 \left| \beta_{n+1} \right| +  \left|
                       \beta_{n-1} \right| \bigg \} \, , \\[1ex]
\\
  \left| \eta_n \right| &\leq&  c \ 
        \bigg \{ \left| \alpha_{n+1} \right| +  \left|
                       \alpha_{n-1} \right| +
                 \left| \beta_{n+1} \right| +  \left|
                       \beta_{n-1} \right| \bigg \} \, ,
\end{eqnarray*}
where $c= \frac{1}{2} \rho_* \left| \Gamma \right| \left| \hat{k}_1 \ p_2
- \hat{k}_2 \ p_1 \right| \ $.  Thus
\begin{eqnarray*}
  \bigg \| \left(
    \begin{array}{cc}
      \xi \\
      \eta
    \end{array}
 \right) \bigg \|^2 \leq 8 c^2 \ \bigg \| 
 \left( 
   \begin{array}{cc}
     \alpha \\
     \beta
   \end{array}
 \right) \bigg \|^2 \, ,
\end{eqnarray*}
which proves the lemma. $\Box$

A complex number $\lambda$ is an eigenvalue of $\LL$, if there
exists $( \alpha, \beta)^T \in \ell_2 \times \ell_2$ such that
\begin{equation}
\LL \left( 
  \begin{array}{cc}
    \alpha \\
    \beta
  \end{array}
 \right) =  \lambda 
 \left(  
   \begin{array}{cc}
     \alpha \\
     \beta
   \end{array}
 \right) \, .
\label{evect}
\end{equation}
\\
\begin{lemma} The eigenvalues of the linear operator $\LL$ have the following
properties: 
\begin{enumerate}
  \item If $\lambda \in C$ is an eigenvalue of $\LL$, then both 
    $\bar{\lambda}$ (the complex conjugate of $\lambda$) and $-\la$ are also
    eigenvalues.

    \item If $\lambda$ is a real eigenvalue, then $\lambda$ is a
      multiple eigenvalue.

    \item  If $\lambda$ is a simple eigenvalue which is not real,
      then its corresponding eigenvector satisfies the relation
      $\beta = \pm \ i \alpha$.
\end{enumerate}
\label{spele1}
\end{lemma}

Proof: Since $\LL$ is a real linear operator, if $\la$ is an eigenvalue 
of $\LL$, then $\bar{\la}$ is also an eigenvalue. Next we show that 
if $\la$ is an eigenvalue 
of $\LL$, then $-\la$ is also an eigenvalue. Let $\la$ and $\om$ be an 
eigenvalue and a corresponding eigenvector, starting from the equation 
(\ref{CLE}), we have
\\
\begin{eqnarray}
\la \omega_{\hat{k} + np} &=& A(p, \hat{k} + (n-1) p) 
     \ \Gamma \ \omega_{\hat{k} + (n-1) p} \nonumber \\ 
& & + \ A(-p, \hat{k} + (n+1)p)\ 
     \bar{\Gamma} \ \omega_{\hat{k} +(n+1)p}\ . \nonumber
\end{eqnarray}
\\
\nid
Then,
\[
\hat{\om}_{\hat{k} + np} = (-1)^n \omega_{\hat{k} + np}
\]
satisfies
\\
\begin{eqnarray}
(-\la) \hat{\omega}_{\hat{k} + np} &=& A(p, \hat{k} + (n-1) p) 
     \ \Gamma \ \hat{\omega}_{\hat{k} + (n-1) p} \nonumber \\ 
& & + \ A(-p, \hat{k} + (n+1)p)\ 
     \bar{\Gamma} \ \hat{\omega}_{\hat{k} +(n+1)p}\ . \nonumber
\end{eqnarray}
\\
\nid
Therefore, $-\la$ is also an eigenvalue. 
To prove claims~2 and~3, notice that the
  Hamiltonian $\HH_{\hat{k}}$ (\ref{LHR}) is invariant under the
  transformation
\\
  \begin{eqnarray*}
\left \{
    \begin{array}{rcl}
\tilde{\alpha}_n &=& - \beta_n \, , \\ \\ 
\tilde{\beta}_n &=& \alpha_n \, .
    \end{array}
\right.
  \end{eqnarray*}
\\
\nid
Therefore, if $\lambda$ is an eigenvalue and
\\
\begin{eqnarray*}
  e^{\lambda t} \left( 
    \begin{array}{cc}
      \alpha \\
      \beta
    \end{array}
 \right)
\end{eqnarray*}
\\
\nid
solves the linear system (\ref{LHR}), then 
\\
\begin{eqnarray*}
    e^{\lambda t} \left( 
      \begin{array}{cc}
        - \beta \\
        \alpha 
      \end{array}
\right)
\end{eqnarray*}
\\
\nid
also solves the system (\ref{LHR}).  If $\lambda$ is real, then
$(\alpha, \beta)^T$ can be chosen to be real.  Assume that
\\
\begin{eqnarray}
   \left(
    \begin{array}{cc}
      \alpha \\
      \beta
    \end{array}
\right) = \eta \left( 
  \begin{array}{cc}
    - \beta \\
    \alpha
  \end{array}
 \right) \, ,
\label{pfe1}
\end{eqnarray}
\\
\nid
then $\eta^2=-1$ which leads to a contradiction.  Thus $( \alpha ,
\beta)^T$ and $(- \beta , \alpha)^T$ are linearly independent,
and $\lambda$ is a multiple eigenvalue which proves claim~2.  If
$\lambda$ is simple and not real, then $(\alpha , \beta)^T$
and $(- \beta  , \alpha)^T$ are linearly dependent and
(\ref{pfe1}) implies that $\beta = \pm \ i \alpha$, which proves claim~3.
$\Box$

In fact, we have the following theorem.
\begin{theorem}
  The eigenvalues of the linear operator $\LL$ defined in (\ref{spl}), have the
  following properties:
  \begin{enumerate}
  \item The eigenvalues of $\LL$ are of 
four types: real pairs ($c, -c$), purely imaginary pairs ($id, -id$), 
quadruples ($\pm c \pm id$), and zero eigenvalues.
    \item If $\lambda$ is a real eigenvalue, then $\lambda$ is
        a multiple eigenvalue.
      If $\lambda$ is a simple non-real eigenvalue, then
        its corresponding eigenvector satisfies the relation
        $\beta = \pm \ i \alpha$.
   \item  If $\sum_{\hat{k}} \cap \bar{D}_{\left| p \right|} =
    \emptyset$, then all the eigenvalues of $\LL$ are either purely
    imaginary and in complex conjugate pairs ($\lambda = ic ,
    -ic$; $c$ is real and not zero) or zeros.
\item  If $| \hat{k} | = | p |$, then
    zero is a multiple eigenvalue of $\LL$.
\item  If $| \hat{k}| =| p |$ and
    $\hat{k} + p \notin \bar{D}_{\left| p \right|}$ (or $\hat{k}
    -p \notin \bar{D}_{\left| p \right|}$), then all the
    eigenvalues for the system (\ref{CLE1}) (or (\ref{CLE2})) are
    either purely imaginary and in complex conjugate pairs or zeros.
  \end{enumerate}
\label{prpev}
\end{theorem}

Proof:
  Claims 1 and 2 are proved in Lemma~\ref{spele1}. Next we prove claim~3.  
If $\sum_{\hat{k}}
  \cap \bar{D}_{\left| p \right|} = \emptyset$, then $I_{\hat{k}}$
  is negative definite for $(\alpha, \beta)^T \in \ell_2
  \times \ell_2$.  In fact, there exists a constant $\delta >0$,
  such that
\\
\begin{equation}
    \delta < - \rho_n <2 \, , \quad  \forall n \in Z \, .
\label{spf0}
\end{equation}
\\
\nid
Thus,
\\
\begin{equation}
  -I_{\hat{k}} > \delta \ \| ( \alpha , \beta )^T \|^2 \, .
\label{spf1}
\end{equation}
\\
\nid
Let $\lambda$ be an eigenvalue of $\LL$, and $(\tilde{\alpha},
\tilde{\beta})^T$ be its corresponding eigenvector, which are written in
terms of real and imaginary parts,
\\
\begin{equation}
   \lambda = \lambda_r + i \ \lambda_i \, , \quad
  (\tilde{\alpha}, \tilde{\beta})^T 
     = (\tilde{\alpha}^{(1)}, \tilde{\beta}^{(1)})^T 
       + i\ (\tilde{\alpha}^{(2)}, \tilde{\beta}^{(2)})^T \, . 
\label{spf2}
\end{equation}
\\
\nid
Then both the real and the imaginary parts of $e^{\lambda
  t}(\tilde{\alpha} , \tilde{\beta})^T$, denoted by $(
\alpha^{(1)} , \beta^{(1)})$ and $( \alpha^{(2)} , \beta^{(2)})$,
are real solutions to system (\ref{LHR}), where
\\
\begin{eqnarray*}
  \alpha^{(1)} &=& e^{\lambda_rt}
                   \left[ \tilde{\alpha}^{(1)} \cos \lambda_it -
                           \tilde{\alpha}^{(2)} \sin \lambda_it
                   \right] \, , \\[1ex]
\\
  \beta^{(1)} &=& e^{\lambda_rt}
                   \left[ \tilde{\beta}^{(1)} \cos \lambda_it -
                           \tilde{\beta}^{(2)} \sin \lambda_it
                   \right] \, ; \\[1ex]         
\\
  \alpha^{(2)} &=& e^{\lambda_rt}
                   \left[ \tilde{\alpha}^{(1)} \sin \lambda_it +
                           \tilde{\alpha}^{(2)} \cos \lambda_it
                   \right] \, , \\[1ex]
\\
  \beta^{(2)} &=& e^{\lambda_rt}
                   \left[ \tilde{\beta}^{(1)} \sin \lambda_it +
                           \tilde{\beta}^{(2)} \cos \lambda_it
                   \right] \, .
\end{eqnarray*}
\\
\nid
Denote by $I^{(j)}_{\hat{k}}$ $ (j=1,2)$ the values of
$I_{\hat{k}}$ evaluated at $(\alpha^{(j)} , \beta^{(j)})^T$;
then
\\
\begin{eqnarray*}
  I^{(1)}_{\hat{k}}(t) + I^{(2)}_{\hat{k}}(t)
     &=& e^{2 \lambda_r t} \sum_{n \in Z} \rho_n
           \left[ (\tilde{\alpha}^{(1)})^2 + (\tilde{\alpha}^{(2)})^2 +
                   (\tilde{\beta}^{(1)})^2 + (\tilde{\beta}^{(2)})^2 
            \right] \\[1ex]
\\
     &=&   I^{(1)}_{\hat{k}}(0) + I^{(2)}_{\hat{k}}(0) \\[1ex]
\\
     &=&   \sum_{n \in Z} \rho_n 
            \left[ (\tilde{\alpha}^{(1)})^2 + (\tilde{\alpha}^{(2)})^2 +
                   (\tilde{\beta}^{(1)})^2 + (\tilde{\beta}^{(2)})^2 
            \right] \\[1ex]
\\
     & \ne& 0 \, ,
\end{eqnarray*}
\\
\nid
by Lemma~\ref{const} and relations (\ref{spf0}, \ref{spf1}).
Thus,
\\
\begin{displaymath}
  e^{2 \lambda_rt} = 1 \, , \quad \hbox{i.e.\ } \lambda_r = 0 \, .
\end{displaymath}
\\
\nid
This proves claim~3.  To prove claim~4, notice that if 
$|\hat{k}| = |p|$, then system (\ref{CLE})
decomposes into systems (\ref{CLE1}, \ref{CLE2}, \ref{CLE3}).
Thus the vectors $(\alpha , \beta)^T$ defined as
\\
\begin{displaymath}
  \alpha_0 = 1 \, , \qquad  
  \alpha_n = 0 \quad (n \ne 0) \, , \qquad
  \beta_n = 0 \quad (\forall n \in Z) \, ;
\end{displaymath}
\\
\nid
and
\\
\begin{displaymath}
  \beta_0 = 1 \, , \qquad 
  \beta_n = 0  \quad (n \ne 0) \, , \qquad
  \alpha_n = 0 \quad  (\forall n \in Z) \, ;
\end{displaymath}
\\
\nid
give two linearly independent eigenvectors of $\LL$ with
eigenvalue zero.  This proves claim~4.
Claim~5 follows from the proof for claim~3 when restricted
to system (\ref{CLE1}) or (\ref{CLE2}). $\Box$
\begin{remark}
  For a finite dimensional linear Hamiltonian system, it is well-known
that the eigenvalues are of 
four types: real pairs ($c, -c$), purely imaginary pairs ($id, -id$), 
quadruples ($\pm c \pm id$), and zero eigenvalues \cite{Poi99} \cite{Lia49} 
\cite{Arn80}. There is also a complete theorem on the normal forms of such
Hamiltonians \cite{Arn80}.
\end{remark}

\newpage
\eqnsection{The Point Spectrum of the Linearized Two-Dimensional Euler 
Equation: A Continued Fraction Study}

Rewrite the equation (\ref{CLE}) as follows,
\begin{equation}
\rho^{-1}_n \dot{\tz}_n = a\ \bigg [ \tz_{n+1} - \tz_{n-1} \bigg ]\ ,
\label{cfr1}
\end{equation}
where $\tz_n = \rho_n e^{in (\th +\pi /2)} \om_{\hat{k}+np}$, $\th +\ga 
=\pi /2$, $\Ga = |\Ga| e^{i\ga}$, $a = {1 \over 2} |\Ga| \left | \begin{array}
{lr} p_1 & \hat{k}_1 \\ p_2 & \hat{k}_2 \\ \end{array} \right |$, 
$\rho_n = | \hat{k}+np|^{-2} - |p|^{-2}$. Let $\tz_n = e^{\la t} z_n$, where 
$\la \in C$; then $z_n$ satisfies 
\begin{equation}
a_n z_n +z_{n-1} - z_{n+1} = 0 \ , 
\label{cfr2}
\end{equation}
where $a_n = \la (a \rho_n)^{-1}$. Let $w_n = z_n / z_{n-1}$ \cite{MS61}; 
then $w_n$ satisfies
\begin{equation}
a_n + {1 \over w_n} = w_{n+1}\ .
\label{cfr3}
\end{equation}
Iteration of (\ref{cfr3}) leads to the continued fraction solution \cite{MS61},
\begin{equation}
w_n^{(1)}=a_{n-1} +{1 \over a_{n-2} + {1 \over a_{n-3}+_{\ \ddots}}}\ \ .
\label{cfr4}
\end{equation}
Rewrite (\ref{cfr3}) as follows,
\begin{equation}
w_n = {1 \over -a_n +w_{n+1}}\ . 
\label{cfr5}
\end{equation}
Iteration of (\ref{cfr5}) leads to the continued fraction solution 
\cite{MS61},
\begin{equation}
w_n^{(2)}=-{1 \over a_{n} + {1 \over a_{n+1}+{1 \over a_{n+2}
+_{\ \ddots}}}}\ \ .
\label{cfr6}
\end{equation}
Before we study the two continued fraction solutions (\ref{cfr4}) and 
(\ref{cfr6}), we like to quote two theorems on the convergence of continued 
fractions \cite{LW92}.
\begin{theorem} [\'{S}leszy\'{n}ski-Pringsheim's Theorem]
The continued fraction
\[
{\bf K}(a_n/b_n)={a_1 \over b_1 + {a_2 \over b_2+{a_3 \over b_3
+_{\ \ddots}}}}\ \ ,
\]
where $\{ a_n \}$ and $\{ b_n \}$ are complex numbers and all $a_n \neq 0$, 
converges if for all $n$
\[
|b_n| \geq |a_n| +1\ .
\]
Under the same condition
\[
|{\bf K}(a_n/b_n)| \leq 1 \ .
\]
\end{theorem}
\begin{theorem} [Van Vleck's Theorem]
Let $0 < \e < \pi/2$, and let $b_n$ satisfy
\[
-\pi/2 +\e < \arg \{ b_n \} < \pi/2 -\e
\]
for all $n$. Then the continued fraction ${\bf K}(1/b_n)$ converges if 
and only if 
\[
\sum_{n=1}^{\infty} |b_n| = \infty \ .
\]
\end{theorem}
Notice that as $n \ra \pm \infty$,
\begin{equation}
a_n \ra \ta = - \la a^{-1} |p|^2 \ . \label{cfr7}
\end{equation}
Then we have the corollary.
\begin{corollary}
If $\mbox{Re}\{ \ta \} \neq 0$, or $\mbox{Re}\{ \ta \}=0 \ (|\ta| >2)$, 
then the two continued fractions (\ref{cfr4}) and (\ref{cfr6}) converge.
\label{cgtle}
\end{corollary}
Proof: If $\mbox{Re}\{ \ta \} \neq 0$, then there exists an positive 
integer $\tN$ and a positive constant $\e$, such that
\[
-\pi /2 +\e < \arg \{ a_n \} < \pi /2 - \e
\]
for all $|n| \geq \tN$, or 
\[
-\pi /2 +\e < \arg \{ -a_n \} < \pi /2 - \e
\]
for all $|n| \geq \tN$. In either case, applying Van Vleck's theorem, we 
have the convergence of the two continued fractions (\ref{cfr4}) and 
(\ref{cfr6}). If $\mbox{Re}\{ \ta \}=0 \ (|\ta| >2)$, then there exists 
an positive integer $\hat{N}$ such that 
\[
|a_n| >2
\]
for all $|n| \geq \hat{N}$. Then applying \'{S}leszy\'{n}ski-Pringsheim's 
theorem, we have the convergence of the two continued fractions (\ref{cfr4}) 
and (\ref{cfr6}). $\Box$
\begin{remark}
In fact, as proved in Theorem \ref{spthla}, $\mbox{Re}\{ \ta \}=0$ and 
$|\ta| \leq 2$ correspond to the continuous spectrum (= essential spectrum) of 
the system.
\end{remark}

When $\mbox{Re}\{ \ta \} \neq 0$, or $\mbox{Re}\{ \ta \}=0 \ (|\ta| >2)$,
as $n \ra - \infty$,
\begin{equation}
w_n^{(1)} \ra w^{(1)}=\ta + {1 \over \ta + {1 \over \ta +_{\ \ddots}}}
=\ta + {\bf K}(1/\ta)\ ,
\label{cfr8}
\end{equation}
as $n \ra +\infty$,
\begin{equation}
w_n^{(2)} \ra w^{(2)}=-{1 \over \ta + {1 \over \ta +{1 \over \ta
+_{\ \ddots}}}}
=-{\bf K}(1/\ta)\ .
\label{cfr9}
\end{equation}
Both $w^{(1)}$ and $w^{(2)}$ satisfy the equation,
\begin{equation}
\tw^2 - \ta \tw -1 = 0\ .
\label{cfr10}
\end{equation}
\begin{itemize}
\item When $\mbox{Re}\{ \ta \} \neq 0$, the solutions of (\ref{cfr10}) 
can be written as:
\begin{equation}
w_{\pm} = {1 \over 2}\bigg [\ta \pm \dl \sqrt{\ta^2 +4}\bigg ]\ ,
\label{cfr11}
\end{equation}
where $\dl = \mbox{sign}(\mbox{Re} \{ \ta \}) \ \mbox{sign}(\mbox{Re} 
\{ \sqrt{\ta^2+4} \})$. (Note that if $\mbox{Re} \{ \ta \} \neq 0$, 
then $\mbox{Re} \{ \sqrt{\ta^2+4} \} \neq 0$.) 
\item When $\mbox{Re}\{ \ta \}=0 \ (|\ta| >2)$, let $\ta = i \xi$, $\xi$ 
is a real number, the solutions of (\ref{cfr10}) can be written as:
\begin{equation}
w_{\pm} = {i \over 2} \bigg [ \xi \pm \dl \sqrt{\xi^2 -4}\bigg ]\ ,
\label{cfr12}
\end{equation}
where $\dl = \mbox{sign} (\xi) \ \mbox{sign}(\sqrt{\xi^2 -4})$.
\end{itemize}
\begin{lemma}
The solutions (\ref{cfr11}) and (\ref{cfr12}) satisfy the inequality
\[
|w_+| > 1 > |w_-|,
\]
and the continued fractions (\ref{cfr8}) and (\ref{cfr9}) have the 
values
\[
w^{(1)} = w_+\ , \ \ w^{(2)} = w_-\ .
\]
\label{cfrle1}
\end{lemma}
Proof: First we show that $|w_+| > 1$ when $\mbox{Re}\{ \ta \} \neq 0$.
\begin{equation}
w_+ \overline{w_+} = {1 \over 4} \bigg ( |\ta|^2 + |\ta^2 +4| + 
\bar{\ta} \dl \sqrt{\ta^2 +4} + \ta \dl \sqrt{\bar{\ta}^2 +4} \bigg ) \ .
\label{cfr13}
\end{equation}
Let $\ta = a_1 + i a_2$, $\sqrt{\ta^2 +4} = b_1 +i b_2$; then
\begin{eqnarray}
b_1 b_2 &=& a_1 a_2 \ , \label{cfr14} \\ 
b_1^2 - b_2^2 &=& a_1^2 -a_2^2 +4 \ . \label{cfr15} 
\end{eqnarray}
From (\ref{cfr14}), we have 
\[
a_1^2 (a_2b_2) = b_2^2 (a_1b_1) \ .
\]
Thus, $a_2b_2$ is either zero or of the same sign as $a_1b_1$. Therefore,
\begin{equation}
\bar{\ta} \dl \sqrt{\ta^2 +4} + \ta \dl \sqrt{\bar{\ta}^2 +4}
=2\dl [a_1b_1+a_2b_2] > 0 \ .
\label{cfr16}
\end{equation}
Together with 
\[
|\ta^2 +4| \geq 4-|\ta|^2\ ,
\]
we have $|w_+| > 1$. When $\mbox{Re}\{ \ta \}=0 \ (|\ta| >2)$, it is obvious 
that $|w_+| > 1$. Notice that $w_+ w_- = -1$, we have $|w_-| < 1$ in 
both cases. Thus, we have
\[
|w_+| > 1 > |w_-| \ .
\]
From the relation $w_+ w_- = -1$, when $\mbox{Re}\{ \ta \} \neq 0$, $
\mbox{Re} \{ w_+ \}$ and $ \mbox{Re} \{ w_- \}$ are of opposite signs.   
When $\mbox{Re}\{ \ta \} \neq 0$, $\mbox{Re} \{ w^{(1)} \}$ and 
$ \mbox{Re} \{ w^{(2)} \}$ are of opposite signs, and $\mbox{Re} 
\{ w^{(1)} \}$ and $ \mbox{Re} \{ w_+ \}$ are of the same sign; thus, 
$w^{(1)} = w_+$ and $w^{(2)} = w_-$. When $\mbox{Re}\{ \ta \}=0 \ (|\ta| >2)$, 
by \'{S}leszy\'{n}ski-Pringsheim's theorem,
\[
|{\bf K}(1/\ta)| \leq 1 \ ;
\]
then,
\[
|w^{(1)}| = |\ta + {\bf K}(1/\ta)| \geq |\ta| - |{\bf K}(1/\ta)| > 1 \ .
\]
Thus, $w^{(1)} = w_+$ and $w^{(2)} = w_-$. $\Box$
\begin{definition}
Define $w^{(*)}$ as follows
\[
w^{(*)}_n = w_n^{(1)}\ \ \mbox{for}\ n \leq 1, \ \ \ 
w^{(*)}_n = w_n^{(2)}\ \ \mbox{for}\ n > 1.
\]
\end{definition}
Then $w^{(*)}_n$ solves (\ref{cfr3}), provided that $w_1^{(1)}=w_1^{(2)}$, 
i.e. 
\begin{equation}
f= a_0 + \bigg ( {1 \over a_{-1} + {1 \over a_{-2} +{1 \over a_{-3}
+_{\ \ddots}}}} \bigg ) + \bigg ( {1 \over a_{1} + {1 \over a_{2} +
{1 \over a_{3}+_{\ \ddots}}}} \bigg ) = 0 \ ,
\label{cfr17}
\end{equation}
where $f = f(\tla,\hat{k},p)$, $\tla = \la /a$. Let $z^{(*)}$ satisfy
\[
w^{(*)}_n = z^{(*)}_n / z^{(*)}_{n-1};
\]
then as $n \ra + \infty$,
\[
z^{(*)}_n \sim (w_-)^n \ ,
\]
and as $n \ra - \infty$,
\[
z^{(*)}_n \sim (w_+)^n \ .
\]
Thus by Lemma \ref{cfrle1}, $z^{(*)} \in \ell_2$. Therefore equation 
(\ref{cfr17}) determines eigenvalues.

If $\Sg_{\hat{k}} \cap \bar{D}_{|p|} = \emptyset$, then $\rho_n < 0$ 
for any $n \in Z$. If $\mbox{Re}\{ \tla \} \neq 0$, then $\mbox{Re}
\{ \ta_n \} \neq 0$ and are of a fixed sign for any $n \in Z$. Then 
$\mbox{Re}\{ f \} \neq 0$. Therefore, in such cases, there is no 
eigenvalue with nonzero real part. This fact is already obtained in Theorem 
\ref{prpev}. Moreover, we have the following fact.
\begin{lemma}
If $\Sg_{\hat{k}} \cap \bar{D}_{|p|} = \emptyset$, then equation (\ref{cfr17}) 
determines no eigenvalue.
\label{lenev}
\end{lemma}
Proof: As discussed above, if $\Sg_{\hat{k}} \cap \bar{D}_{|p|} = \emptyset$, 
then the possible solution $\tla$ to (\ref{cfr17}) has to be imaginary. 
Therefore, $\ta$ has to be imaginary. Rewrite equation (\ref{cfr2}) as 
follows
\begin{equation}
\tilde{L} z_n \equiv {\rho_n \over \rho} [ z_{n+1} -z_{n-1}] = \ta z_n \ .
\label{cfr18}
\end{equation}
If $\Sg_{\hat{k}} \cap \bar{D}_{|p|} = \emptyset$, then 
$0 < \rho_n / \rho < 1$ for all $n \in Z$. Thus, $\| \tilde{L} \| \leq 2$. 
Then if $\tla$ is an eigenvalue, $|\ta| \leq \| \tilde{L} \| \leq 2$. 
Notice that equation (\ref{cfr17}) should be solved under the condition 
$\mbox{Re}\{ \ta \} \neq 0$ or $\mbox{Re}\{ \ta \}=0 \ (|\ta| >2)$; thus, 
in this case, equation (\ref{cfr17}) determines no eigenvalue. $\Box$
\begin{remark}
In fact, as proved in Theorem \ref{spthla}, $\mbox{Re}\{ \ta \}=0$ and 
$|\ta| \leq 2$ correspond to the continuous spectrum (= essential spectrum) of 
the system. Thus, if $\Sg_{\hat{k}} \cap \bar{D}_{|p|} = \emptyset$, the 
point spectrum is empty (cf: Theorem \ref{spthla}).
Then the problem is reduced to solving equation (\ref{cfr17}) under the 
conditions $\Sg_{\hat{k}} \cap \bar{D}_{|p|} \neq \emptyset$, 
$\mbox{Re}\{ \ta \} \neq 0$ or $\mbox{Re}\{ \ta \}=0 \ (|\ta| >2)$.
\end{remark}

\begin{figure}[ht]
  \begin{center}
    \leavevmode
      \setlength{\unitlength}{2ex}
  \begin{picture}(36,27.8)(-18,-12)
    \thicklines
\put(0,-14){\vector(0,1){28}}
\put(-18,0){\vector(1,0){36}}
\put(0,15){\makebox(0,0){$\Im \{ \tla \}$}}
\put(18.5,0){\makebox(0,0)[l]{$\Re \{ \tla \}$}}
\put(2.4,3.5){\circle*{0.5}}
\put(-2.4,3.5){\circle*{0.5}}
\put(2.4,-3.5){\circle*{0.5}}
\put(-2.4,-3.5){\circle*{0.5}}  
\end{picture}
  \end{center}
\caption{The quadruple of eigenvalues determined by equation (\ref{cfr17})
for the system led by the class $\Sg_{\hat{k}}$ labeled by $\hk = (1,0)^T$, 
when $p=(1,1)^T$.}
\label{figev}
\end{figure}
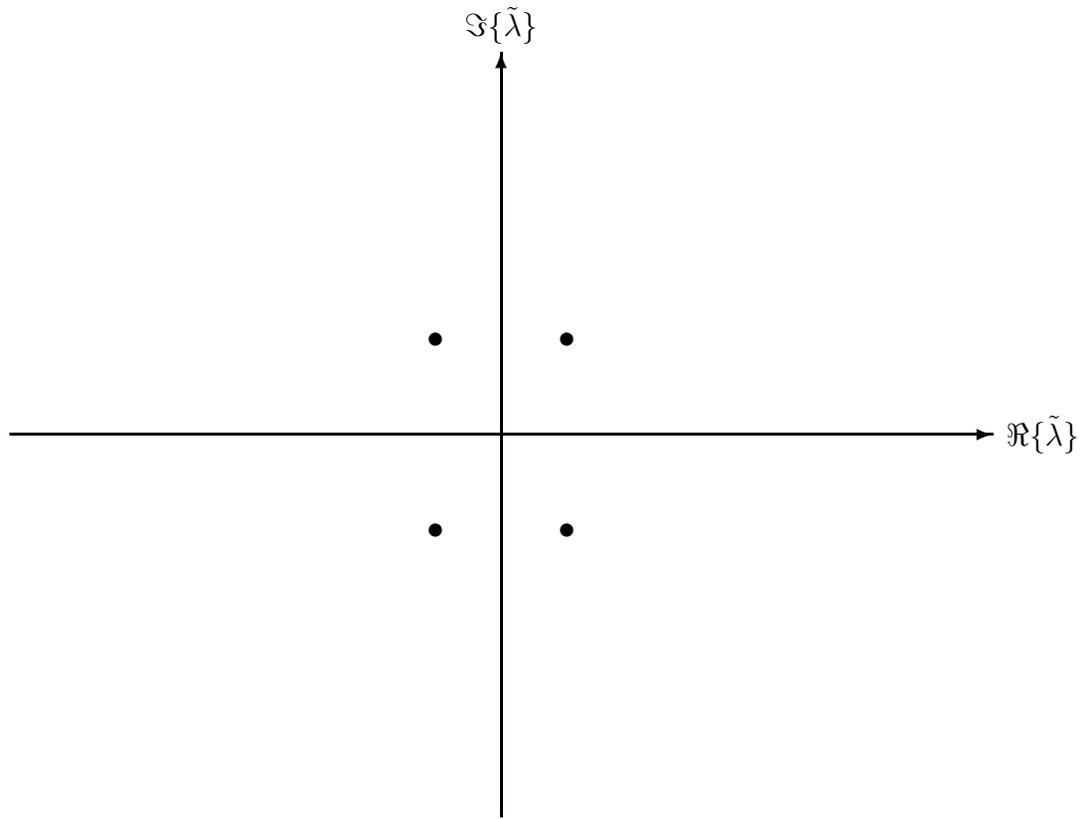

{\bf Example}: Let $p=(1,1)^T$, in this case, only one class $\Sg_{\hat{k}}$ 
labeled by $\hk = (1,0)^T$ has no empty intersection with $\bar{D}_{|p|}$ 
(the other class labeled by $\hk = (0,1)^T$ gives the complex conjugate 
of the system led by the class labeled by $\hk = (1,0)^T$). For this class,
$| \rho_n/\rho | \leq 1$ for all $n \in Z$. Thus, the linear operator 
$\tilde{L}$ defined in (\ref{cfr18}) has norm $\| \tilde{L} \| \leq 2$.
Therefore, equation (\ref{cfr17}) determines no real eigenvalue. Numerical 
calculation on equation (\ref{cfr17}) gives the eigenvalue:
\[
\tla=0.24822302478255 \ + \ i \ 0.35172076526520\ .
\]
By Theorem \ref{prpev}, equation (\ref{cfr17}) determines a quadruple 
of eigenvalues, see figure \ref{figev} for an illustration.

\newpage
\eqnsection{The Spectra of the Linearized Two-Dimensional Euler
Equation: An Infinite Matrix Study}

\subsection{The General Setup}

Rewrite (\ref{CLE}) as follows:
\\
\begin{equation}
  \label{ler1}
  \dot{\tilde{z}}_n = ia \ \bigg [ \rho_{n-1} \ \tilde{z}_{n-1} +
    \rho_{n+1} \ \tilde{z}_{n+1} \bigg ] \, , \quad (n \in Z)
\end{equation}
\\
\nid
where 
\\
\begin{eqnarray*}
  \tilde{z}_n = e^{in \theta} \omega_{\hat{k}+np} \, , \qquad
  \Gamma = \left| \Gamma \right| e^{i\ga} \, , \qquad
  \theta + \gamma = \pi/2 \, , \qquad
  a = \frac{1}{2} \left| \Gamma \right| 
\left|
  \begin{array}{cc}
p_1 & \hat{k}_1 \\
p_2 & \hat{k}_2
  \end{array}
\right| \, .
\end{eqnarray*}
\\
\nid
Relabel $\left\{ \tilde{z}_n \right\}$ as follows:
\\
\begin{eqnarray*}
  \left\{
    \begin{array}{rcl@{\qquad}rcl}
\tilde{z}_n &=& z_{2n} \, ,      & n &\geq& 1 \, , \\
\tilde{z}_{-n} &=& z_{2n+1} \, , & n &\geq& 0 \, ;
    \end{array}
  \right.
\end{eqnarray*}
\\
\nid
then
\\
\begin{eqnarray}
\dot{z}_{2n} &=& ia \ 
        \bigg [ \rho_{n-1} \ z_{2(n-1)} + \rho_{n+1} \ z_{2(n+1)} 
        \bigg ] \, , \quad (n \geq 2) \label{ler2} \\
\nonumber \\
\dot{z}_{2n+1} &=& ia \ 
         \bigg [ \rho_{-n+1} \ z_{2(n-1)+1} + \rho_{-n-1}
           \ z_{2(n+1)+1}
         \bigg ] \, , \quad (n \geq 1) \label{ler3} \\
\nonumber \\
\dot{z}_2 &=& ia \ 
         \bigg [ \rho_0 \ z_1 + \rho_2 \ z_4
         \bigg ] \, , \label{ler4} \\ 
\nonumber \\
\dot{z}_1 &=& ia \ 
          \bigg [ \rho_1 \ z_2 + \rho_{-1} \  z_3
          \bigg ] \, , \label{ler5}
\end{eqnarray}
\\
\nid
for $z=(z_1, z_2, \cdot \cdot \cdot )^T$.  Notice that Equations (\ref{ler2})
and (\ref{ler3}) are decoupled, the coupling between components
of $z$ with even and odd indices is through Equations (\ref{ler4})
and (\ref{ler5}).  The right hand side of 
(\ref{ler2})--(\ref{ler5}) define a bounded linear operator $\LL_A : \ell_2
\mapsto \ell_2$, with the infinite matrix representation, 
\\
\begin{eqnarray}
A= ia \left(
    \begin{array}{cccccccc}
      0 &     \rho_1 & \rho_{-1} & 0 & 0 & 0 & 0 \\
      \rho_0 & 0 & 0 & \rho_2 & 0 & 0 & 0 \\
      \rho_0 & 0 & 0 & 0 & \rho_{-2} & 0 & 0 & \bigcirc \\
      0 & \rho_1 & 0 & 0 & 0 & \rho_3 & 0 \\
      0 & 0 & \rho_{-1} & 0 & 0 & 0 & \rho_{-3} \\
      & \bigcirc & & \ddots && \bigcirc & & \ddots
    \end{array}
       \right) \, .
\label{mata}
\end{eqnarray}
\\
\nid
More importantly,
\\
\begin{equation}
     \rho_n \longrightarrow \rho = - \left| p \right|^{-2} \, ,
     \quad \hbox{as } \left| n \right| \to \infty \, .
\label{lim}
\end{equation}
\\
\nid
Define the infinite matrix
\\
\begin{eqnarray}
B= ib \left(
    \begin{array}{cccccccc}
      0 & 1 & 1 & 0 & 0 & 0 & 0 \\
      1 & 0 & 0 & 1 & 0 & 0 & 0 \\
      1 & 0 & 0 & 0 & 1 & 0 & 0 &\bigcirc \\
      0 & 1 & 0 & 0 & 0 & 1 & 0 \\
      0 & 0 & 1 & 0 & 0 & 0 & 1 \\
      & \bigcirc & & \ddots && \bigcirc & & \ddots
    \end{array}
       \right) \, ,
\label{matb}
\end{eqnarray}
\\
\nid
where $b=a \rho = - a \left| p \right|^{-2}$.
Define the infinite matrix $C$ as
\\
\begin{equation}
  C=A-B \, , \label{matc1}
\end{equation}
\\
\nid
that is,
\\
\begin{eqnarray}
C= ia \left(
    \begin{array}{c@{\quad}c@{\quad}c@{\quad}c@{\quad}c@{\quad}c@{\quad}cc}
      0 & \tilde{\rho}_1 & \tilde{\rho}_{-1}  & 0 & 0 & 0 & 0&  \\
      \tilde{\rho}_0  & 0 & 0 & \tilde{\rho}_2 & 0 & 0 & 0& \\
      \tilde{\rho}_0  & 0 & 0 & 0 & \tilde{\rho}_{-2}  &0 &0 & \bigcirc \\
      0 & \tilde{\rho}_1 & 0 & 0 & 0 & \tilde{\rho}_3 & 0 & \\
      0 & 0 & \tilde{\rho}_{-1} & 0 & 0 & 0 & \tilde{\rho}_{-3} &\\
      & \bigcirc & &\ddots& & \bigcirc & &\ddots \\
    \end{array}
       \right) \, .
\label{matc2}
\end{eqnarray}
\\
\nid
where $\tilde{\rho}_n = \rho_n -\rho$.
Denote by $\LL_B$ and $\LL_C$ the bounded linear operators with
the infinite matrix representations by $B$ and $C$. According to
Duren \cite{Dur60}, $\LL_A$, $\LL_B$ and $\LL_C$ are called $( 2
\times 2+1)$-operators, since their entries $c_{n,n+m}$ satisfy
the condition $c_{n, n+m} =0$ if $\left| m \right| >2$.  $\LL_B$
is a $(2 \times 2+1)$-operator with constant coefficients, since
its entries $c_{n, n+m}$ is independent of $n$ when $n>2$; and $i\LL_B$
is self-adjoint.

\begin{theorem}
The bounded linear operator $\LL_C : \ell_2 \mapsto \ell_2$ is a
compact operator.
\label{cpt}
\end{theorem}

Proof:
  Denote by $\LL^{(N)}_C$ the linear operator represented through
  the matrix $C_{N \times N}$ obtained from $C$ by replacing its
  entries $c_{m,n}$ by $0$, when $m>N$.  Let $\left\{ z^{(j)}
  \right\}$ be a bounded sequence in $\ell_2$; then $\left\{
    \LL^{(1)}_C z^{(j)} \right\}$ is a bounded sequence in which
  each element has only one nonzero component, i.e.
\\
  \begin{displaymath}
    ( \LL^{(1)}_C z^{(j)} )_n = 0 \, , \quad
    \hbox{when } n>1 \, .
  \end{displaymath}
\\
\nid
Thus, there exists a subsequence $\{ z^{(1_j)}\}$,
such that $\{ \LL^{(1)}_C z^{(1_j)} \}$ converges in
$\ell_2$.  Similarly, we can get a subsequence of $\{
z^{(1_j)}\}$, denoted as  $\{ z^{(2_j)} \}$,
such that  $\{ \LL^{(2)}_C z^{(2_j)} \}$ converges in
$\ell_2$, and so on.  Therefore, we have a nested list of
subsequences:
\\
\begin{eqnarray*}
  \begin{array}{ccc}
z^{(1_1)} & z^{(1_2)} & \cdots \cdots  \cdots\\
z^{(2_1)} & z^{(2_2)} & \cdots \cdots \cdots \\
\vdots & \vdots &  \vdots \\ 
\vdots & \vdots &  \vdots \\
\vdots & \vdots &  \vdots \\ 
\end{array}
\end{eqnarray*}
\\
\nid
We choose the subsequence $\{ z^{(n_n)} \}$ of $\{
z^{(j)} \}$, which is the diagonal of the above list.
There exist constants $\zeta$ and $N_0$, such that
\\
\begin{equation}
  \left\| \LL^{(\hat{n})}_C z^{(n_n)} - \LL_C \ z^{(n_n)} \right\| 
  \leq \ \frac{\zeta}{\hat{n}^2} \, , \quad
  \hbox{for all } \hat{n} > N_0 \hbox{ and all } n \, .
\label{cp1}
\end{equation}
\\
\nid
For any $\epsilon >0$, choose $\hat{N}$
large enough, such that 
\\
\begin{equation}
  \frac{\zeta}{\hat{N}^2} < \frac{1}{3} \ \epsilon \, .
\label{cp2}
\end{equation}
\\
\nid
Since the subsequence  $\left\{ \LL^{(\hat{N})}_C
  z^{(\hat{N}_j)} \right\}$ converges, there exists $\tilde{N}$,
such that 
\\
\begin{equation}
   \left\| \LL^{(\hat{N})}_C z^{(\hat{N}_{j_1})} - 
      \LL^{(\hat{N})}_C z^{(\hat{N}_{j_2})} \right\| < \frac{1}{3} \ 
      \epsilon \, , \quad
      \forall j_1 \, , j_2 > \tilde{N}  \, .
\label{cp3}
\end{equation}
\\
\nid
Let $N_1 = \max \left\{ \hat{N} , \tilde{N} \right\}$; then
\\
\begin{equation}
   \left\| \LL^{(\hat{N})}_C z^{(n_n)} - 
      \LL^{(\hat{N})}_C z^{(\tilde{n}_{\tilde{n}})} \right\| < 
      \frac{1}{3} \ \epsilon \, \quad
      \forall n \, , \tilde{n} > N_1  \, .
\label{cp4}
\end{equation}
\\
\nid
Thus
\\
\begin{eqnarray*}
   \left\| \LL_C \ z^{(n_n)} 
          - \LL_C \ z^{(\tilde{n}_{\tilde{n}})} \right\| 
    &\leq& \left\| \LL_C \ z^{(n_n)} - \LL^{(\hat{N})}_C z^{(n_n)} \right \|  \\[1ex]
\\
    & & + \left\| \LL^{(\hat{N})}_C z^{(n_n)} 
             - \LL^{(\hat{N})}_C z^{(\tilde{n}_{\tilde{n}})}
           \right\| \\ \\ 
         & &+ \left\| \LL^{(\hat{N})}_C z^{(\tilde{n}_{\tilde{n}})}
           - \LL_C \ z^{(\tilde{n}_{\tilde{n}})} \right\|  \\[1ex]
\\
    &<& \frac{1}{3}\ \epsilon + \frac{1}{3}\ \epsilon + \frac{1}{3}
      \ \epsilon = \epsilon \, , \quad \forall n \, , \tilde{n} > N_1 \, .
\end{eqnarray*}
\\
\nid
Therefore, $\left\{ \LL_C \ z^{(n_n)} \right\}$ is a Cauchy sequence
in $\ell_2$; thus converges.  This proves that $\LL_C$ is a
compact operator. $\Box$

\begin{remark}
  In fact, a theorem of Achieser and
  Glasmann \cite{AG58} \cite{Dur60} states that a $(2M+1)$-operator
  is compact if and only if its diagonal sequence entries tend to
  zeros, i.e.\ $c_{n, n+m} \to 0$, as $n \to \infty$ for each
  fixed $m, \ \left| m \right| \leq M$.  Here we give the proof for
  self-containedness.
\end{remark}

\subsection{The Spectra of the Linear Operator $\LL_B$}

Next we will follow a theory of Duren \cite{Dur60} to study the
spectra of the constant-coefficient infinite-matrix bounded self-adjoint 
operator $i\LL_B$.

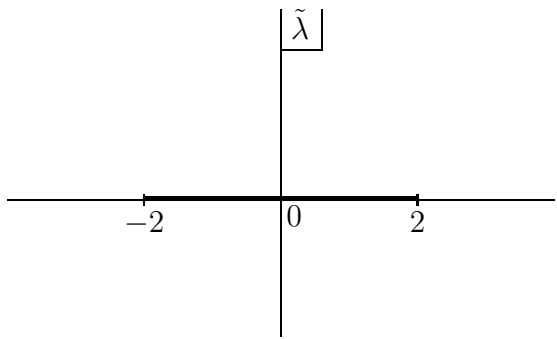
\begin{figure}[ht]
\begin{center}
  \setlength{\unitlength}{2ex}
\begin{picture}(20,12)(-10,-5)
  \setlength{\unitlength}{2ex}
    \thinlines
\put(0,-5){\line(0,1){12}}
\put(-10,0){\line(1,0){20}}
    \thicklines
\put(-5,0.1){\line(1,0){10}}
    \thinlines
\put(-5,-.2){\line(0,1){.4}}
\put(5,-.2){\line(0,1){.4}}
\put(-5,-.4){\makebox(0,0)[t]{$-2$}}
\put(5,-.4){\makebox(0,0)[t]{$2$}}
\put(.2,-.2){\makebox(0,0)[tl]{$0$}}
\put(.7,7){\makebox(0,0)[t]{$\tilde{\lambda}$}}
\put(0,5.5){\line(1,0){1.5}}
\put(1.5,5.5){\line(0,1){1.5}}
\end{picture}
\end{center}
\caption{The spectral curve $C_{\tilde{B}}\ $.}
\label{speccu}
\end{figure}

The {\em{characteristic polynomial}} for the difference equation
\\
\begin{equation}
  (B - \lambda I) z = 0 \, ,
\label{dife}
\end{equation}
\\
\nid
where $I$ is the identity matrix, is defined as:
\\
\begin{equation}
  f_B (w, \lambda) = ib - \lambda\ w^2 + ib\ w^4 \, .
\label{cpl1}
\end{equation}
\\
\nid
Define the rescaled characteristic polynomial as follows:
\\
\begin{equation}
  \tilde{f}_B (w, \tilde{\lambda} ) = 1- \tilde{\lambda}\ w^2 +
  w^4 \, ,
\label{cpl2}
\end{equation}
\\
\nid
where $\lambda = ib \tilde{\lambda}$.  In fact, $\tilde{f}_B (w,
\tilde{\lambda})$ is the characteristic polynomial for the
difference equation
\\
\begin{equation}
  ( \tilde{B} - \tilde{\lambda} I) z=0 \, ,
\label{cpl3}
\end{equation}
\\
where $\tilde{B} =-ib^{-1} B$.  The roots of $\tilde{f}_B (w,
\tilde{\lambda} )$ are:
\\
\begin{equation}
   w_* \, , \quad 
  -w_* \, , \quad
  \frac{1}{w_*} \, , \quad
  - \ \frac{1}{w_*} \, ,
\label{cpl4}
\end{equation}
\\
\nid
where
\\
\begin{displaymath}
  w_* = \left[ 
          \frac{\tilde{\lambda} + \sqrt{\tilde{\lambda}^2 -4}}{2}    
       \right ]^{1/2} \, .
\end{displaymath}
\\
\nid
\begin{definition}  The {\bf spectral curve} of the linear operator
  $\LL_{\tilde{B}}$ (with the infinite matrix 
  representation by $\tilde{B}$), denoted by $C_{\tilde{B}}$, is 
  defined to be
  the set of all $\tilde{\lambda} \in C$ for which the
  characteristic polynomial $\tilde{f}_B (w, \tilde{\lambda})$
  has a root of modulus one.  The {\bf spectral point-set} of
  the operator $\LL_{\tilde{B}}$, denoted by $P_{\tilde{B}}$, is
  defined to be the set of all $\tilde{\lambda} \in C$ for which
  the characteristic polynomial $\tilde{f}_B (w,
  \tilde{\lambda})$ has a multiple root.  Denote by
  $S_{\tilde{B}} ( \tilde{\lambda})$ the number of roots of
  $\tilde{f}_B (w, \tilde{\lambda})$, of modulus less than $1$
  (counted with multiplicity).
\end{definition}
Notice that
\\
\begin{equation}
  \tilde{\lambda} = w^2_* + w^{-2}_* \, .
\label{cpl5}
\end{equation}
\\
\nid
Let $w_*$ be a root of modulus~1 (then all the four roots are of
modulus~1), $w_*=e^{i \theta}$, $\theta \in [0, 2 \pi)$; thus the
spectral curve $C_{\tilde{B}}$ is the segment of the real axis,
\\
\begin{equation}
   C_{\tilde{B}} \ : \ \tilde{\lambda} = 2 \cos 2 \theta \, , \quad
       \theta \in [0, 2 \pi) \, .
\label{cpl6}
\end{equation}
\\
\nid
See Fig.\ref{speccu}.
The spectral point-set $P_{\tilde{B}}$ consists of two points,
\\
\begin{equation}
  P_{\tilde{B}}\ : \ \tilde{\lambda} = \pm \ 2 \, ,
\label{cpl7}
\end{equation}
\\
\nid
which are the boundary points of the spectral curve
$C_{\tilde{B}}$.  At $\tilde{\lambda} = \pm \ 2\ $, the four roots of
$\tilde{f}_B (w, \tilde{\lambda})$ are,
\\
\begin{eqnarray*}
  \begin{array}{lrcl@{\qquad}l}
    \hbox{at } & \tilde{\lambda} &=& 2 :  
        & 1 \, , \ -1 \, , \ 1 \, , \ -1 \, ; \\
    \hbox{at } & \tilde{\lambda} &=& -2 : 
        & i \, , \ -i \, , \ -i \, , \ i \, .
  \end{array}
\end{eqnarray*}
\\
\nid
The function $S_{\tilde{B}}(\tilde{\lambda})$ is
\\
\begin{eqnarray}
  S_{\tilde{B}} (\tilde{\lambda}) = \left\{
    \begin{array}{c@{\quad}l}
      0 \, , & \hbox{if } \ \tilde{\lambda} \in C_{\tilde{B}} \, , \\ \\
      2 \, , & \hbox{if } \ \tilde{\lambda} \notin C_{\tilde{B}} \, . \\
    \end{array}
\right.
\label{cpl8}
\end{eqnarray}
\\
\nid
The general solution to (\ref{cpl3}) is
\\
\begin{eqnarray}
  z_n &=& c_1 w^n_* + c_2(-w_*)^n + c_3 w^{-n}_* + c_4(-w_*)^{-n}
             \, ,  \label{cpl9} \\
      && \qquad \hbox{if } \lambda \notin P_{\tilde{B}} \, , 
         \nonumber \\ \nonumber \\
  z_n &=& c_1 w^n_* + c_2 \ n \ w_*^n + c_3 (-w_*)^n + c_4 \ n\ (-w_*)^{n}
          \, , \label{cpl10}\\
      && \qquad \hbox{if } \lambda \in P_{\tilde{B}} 
            \ (\hbox{then } w_* =1 \, , i) \, ; \nonumber
\end{eqnarray}
\\
\nid
under the restrictions:
\\
\begin{eqnarray}
& &  -\tilde{\lambda} \ z_1 + z_2 + z_3 = 0 \, , \label{cpl11} \\
\nonumber \\
& &  z_1  - \tilde{\lambda} \ z_2 + z_4 = 0 \, .\label{cpl12}
\end{eqnarray}
\\
\begin{lemma}
The general solution (\ref{cpl9}, \ref{cpl10}) is in $\ell_2$
if and only if $\left| w_* \right| <1$ and $c_3 = c_4 =0$ in
(\ref{cpl9}) or $\left| w_* \right| >1$ and $c_1=c_2=0$ in (\ref{cpl9}).
\end{lemma}

Proof:
  See (\cite{Dur60}, pp.~24, Lemma~5). $\Box$

\begin{definition}
Let $\LL : \ell_2 \mapsto \ell_2$ be a linear operator.  The set
of points $\sigma_p(\LL)$ in the complex $\lambda$-plane $C$ such
that $(\LL- \lambda I)$ has no inverse (i.e.\ $\LL - \lambda I$ is
not 1-1), is called the point spectrum of $\LL$.  The set of
points $\sigma_r(\LL)$ in $C$ such that $(\LL - \lambda
I)^{-1}$ exists and is a linear operator with domain not
everywhere dense is called the residual spectrum of $\LL$.  The
set of points $\sigma_c(\LL)$ in $C$ such that $(\LL - \lambda
I)^{-1}$ exists and is an unbounded linear operator with domain
everywhere dense is called the continuous spectrum of $\LL$.  The
set of points $\rho (\LL)$ in $C$ such that $(\LL - \lambda
I)^{-1}$ exists and is a bounded linear operator with domain
everywhere dense is called the resolvent set of $\LL$.  The set
$\sigma (\LL) = \sigma_p(\LL) \cup \sigma_r(\LL) \cup
\sigma_c(\LL)$ is called the spectrum of $\LL$.
\end{definition}

\nid
Without loss of generality, assume $\left| w_* \right| <1$.  Then
$\tilde{\lambda}$ is an eigenvalue if and only if there exists a
non-trivial solution $(c_1, c_2)$ to the following system,
\\
\begin{eqnarray}
  \left(
    \begin{array}{cc}
      - \tilde{\lambda} w_* + w_*^2 + w_*^3  &
      \tilde{\lambda} w_* + w^2_* - w^3_* \\[1ex]
\\
      w_* - \tilde{\lambda} w^2_* + w^4_* &
      -w_* - \tilde{\lambda} w^2_* + w^4_*
    \end{array}
  \right)
  \left(
  \begin{array}{c}
    c_1 \\
    c_2
  \end{array}
   \right) = 0 \, .
\label{cpl13}
\end{eqnarray}
\\
\begin{theorem}
The $\ell_2$ point spectrum $\sigma_p{(B)}$ of the linear operator
$\LL_B$ is empty.
\label{pspec}
\end{theorem}

Proof: The determinant
\\
\begin{eqnarray}
    & &\det \left(
      \begin{array}{cc}
        -\tilde{\lambda} w_* + w^2_* + w^3_*  &
        \tilde{\lambda} w_* + w^2_* - w^3_* \\ \\
        w_* -\tilde{\lambda} w^2_* + w^4_*    &
        -w_*- \tilde{\lambda} w^2_* + w^4_* 
      \end{array}
      \right) \nonumber \\ \nonumber \\ 
     & =& 2 w^3_* \left[ w^4_* - 2 \tilde{\lambda} w^2_* +
  (\tilde{\lambda}^2-1) \right] = 0 \, ,
  \label{pcpl1}
\end{eqnarray}
\\
\nid
implies that
\\
\begin{equation}
  w^4_* - 2 \tilde{\lambda} w^2_* + \tilde{\lambda}^2 -1=0 \, ,
\label{pcpl2}
\end{equation}
\\
\nid
since $w_* \neq 0$.  $w_*$ is a root of $\tilde{f}_B(w,
\tilde{\lambda})$ (\ref{cpl2}),
\\
\begin{equation}
  w^4_* - \tilde{\lambda} w^2_* +1=0 \, .
\label{pcpl3}
\end{equation}
\\
\nid
From (\ref{pcpl2}, \ref{pcpl3}), we have
\\
\begin{equation}
w^2_* = \frac{\tilde{\lambda}^2 -2}{\tilde{\lambda}} \, .
\label{pcpl4}
\end{equation}
\\
\nid
Notice also that
\\
\begin{displaymath}
  \tilde{\lambda} = w^2_* + w^{-2}_* = 
  \frac{\tilde{\lambda}^2 -2}{\tilde{\lambda}} +
  \frac{\tilde{\lambda}}{\tilde{\lambda}^2 -2}\ ,
\end{displaymath}
\\
\nid
which implies that $\tilde{\lambda} = \pm \ 2\ $.  Then $\left| w_*
\right| =1$.  Thus if $\left| w_* \right| <1$, then
Equation (\ref{cpl13}) has only trivial solution.  Therefore, the point
spectrum of $\LL_{\tilde{B}}$ is empty; equivalently, the point
spectrum of $\LL_B$ is empty. $\Box$

\begin{theorem}
The $\ell_2$ residual spectrum $\sigma_r(B)$ of the linear
operator $\LL_B$ is empty.
\label{rspec}
\end{theorem}

Proof:
  $\tilde{\lambda} \in \sigma_r(\tilde{B})$ if and only if the
  dimension of the orthocomplement of $(\LL_{\tilde{B}} -
  \tilde{\lambda}I) \circ \ell_2$ is positive and $(\LL_{\tilde{B}} -
  \tilde{\lambda}I)^{-1}$ exists. From the inner product
  relation
\\
  \begin{eqnarray*}
    \left\langle ( \LL_{\tilde{B}} - \tilde{\lambda}I )
        z^{(1)} , z^{(2)} \right\rangle 
    &=& \left\langle z^{(1)} , \left( \LL^*_{\tilde{B}} -
        \bar{\tilde{\lambda}}I   \right) z^{(2)}  \right\rangle \\ \\ 
    &=& \left\langle z^{(1)} , \left( \LL_{\tilde{B}} -
      \bar{\tilde{\lambda}}I   \right) z^{(2)}  \right\rangle \, ,
  \end{eqnarray*}
\\
\nid
  since $\LL_{\tilde{B}}$ is self-adjoint, $\LL_{\tilde{B}} =
  \LL^*_{\tilde{B}}$ (the adjoint of $\LL_{\tilde{B}}$), where 
$\langle \ , \ \rangle$ denotes the
  inner product over the complex field, we have that if $\tilde{\lambda} \in
  \sigma_r(\tilde{B})$, then $\bar{\tilde{\lambda}}
  \in \sigma_p(\tilde{B})$.  By Theorem~\ref{pspec},
  $\sigma_p(\tilde{B})$ is empty; thus, $\sigma_r(\tilde{B})$
  is empty; equivalently, $\sigma_r(B)$ is empty.
$\Box$ 

Since $i\LL_B$ is self-adjoint, this theorem is well-known, but 
we furnish a short proof here.
From (\ref{dife}, \ref{cpl3}) and (\ref{cpl6}, \ref{cpl7}), the
spectral curve $C_B$ for the linear operator $\LL_B$ is the
segment of the imaginary axis,
\\
\begin{eqnarray}
  C_B \ : \ \lambda = i2b \ \cos 2\theta \, , \quad \theta \in
  \left[0,2 \pi \right) \, ;
\label{srb}
\end{eqnarray}
\\
\nid
the spectral point-set $P_B$ for the linear operator $\LL_B$ is 
\\
\begin{equation}
  P_B \ : \ \lambda = \pm \ i2b \, ,
\label{spb}
\end{equation}
\\
\nid
which are the boundary points of the spectral curve $C_B$.

\begin{theorem}
The $\ell_2$ continuous spectrum $\sigma_c(B)$ of the linear
operator $\LL_B$ is the spectral curve, $\sigma_c(B) = C_B\ $.
The $\ell_2$ resolvent set $\rho(B)$ of the linear operator
$\LL_B$ is the complement of $C_B$ in the finite complex plane
$C$, $\rho(B) = (C_B)'\ $.
 \label{crspec}
\end{theorem}

Proof: First we show that if $\tilde{\lambda} \in C_{\tilde{B}}$, then
  $\tilde{\lambda} \in \sigma_c(\tilde{B})$.
Since both $\sigma_p(\tilde{B})$ and $\sigma_r(\tilde{B})$
are empty by Theorems~\ref{pspec} and \ref{rspec}, for any
$\tilde{\lambda} \in C_{\tilde{B}}$, $(\LL_{\tilde{B}} -
\tilde{\lambda}I)^{-1}$ exists and is everywhere densely
defined.  We need to show that $(\LL_{\tilde{B}} -
\tilde{\lambda}I)^{-1}$ is unbounded.  For any $\tilde{\lambda}
\in C_{\tilde{B}}$, there exists a root of $\tilde{f}_B (w,
\tilde{\lambda})$ of modulus one,
\\
\begin{displaymath}
 w_* = e^{i \theta} \, , \quad \theta \in [0,2 \pi ) \, . 
\end{displaymath}
\\
\nid
Define the elements
\\
\begin{eqnarray*}
  z^{(N)}_n = \left\{
    \begin{array}{cl}
      e^{in \theta} \, , \quad &n \leq N \, , \\ \\
      0 \, , \quad & n>N \, .
    \end{array}
    \right.
\end{eqnarray*}
\\
\nid
Then $z^{(N)} \in \ell_2$ for each finite $N$, and $\| z^{(N)} \|
\to \infty$, as $N \to \infty$.  There exists a constant $d$
independent of $N$, such that
\\
\begin{displaymath}
  \left\| ( \tilde{B} - \tilde{\lambda}I ) z^{(N)} \right\|
  \leq d \, , \quad \forall \ N \, .
\end{displaymath}
\\
\nid
Thus
\\
\begin{displaymath}
  \frac{\left\| z^{(N)} \right\|}{\left\|( \tilde{B} - \tilde{\lambda}I)
       z^{(N)} \right\| } \longrightarrow \infty \, , 
   \quad \hbox{as } N \longrightarrow \infty \, .
\end{displaymath}
\\
\nid
Therefore, $(\LL_{\tilde{B}} - \tilde{\lambda} I)^{-1}$ is
unbounded, and $\tilde{\lambda} \in \sigma_c (\tilde{B})$.  Next
we show that if $\tilde{\lambda} \notin C_{\tilde{B}}$, then
$\tilde{\lambda} \in \rho(\tilde{B})$.  For any
$\tilde{\lambda} \notin C_{\tilde{B}}$, the corresponding roots
of $\tilde{f}_B (w, \tilde{\lambda})$ are (\ref{cpl4}), such that
\\
\begin{equation}
  \left| w_* \right| \ = \ \left| -w_* \right| \ <\ 
  1 \ < \ \left| w^{-1}_* \right|\  = \ \left| (-w_*)^{-1} \right| \, .
\label{pfrv1}
\end{equation}
\\
\nid
For any $ y \in\ell_2$, we want to construct a solution to
\\
\begin{equation}
  (\tilde{B} - \tilde{\lambda}I) z=y \ ,
\label{pfrv2}
\end{equation}
\\
\nid
using the method of variation of coefficients.  Explicitly, we
need to solve 
\\
\begin{equation}
  z_n - \tilde{\lambda} \ z_{n+2} + z_{n+4} = y_{n+2} \, ,
  \quad (n \geq 1)
\label{pfrv3}
\end{equation}
\\
\nid
under the constraints
\\
\begin{eqnarray}
  \left\{
    \begin{array}{rcl}
      - \tilde{\lambda} z_1 + z_2 + z_3 &=& y_1 \, , \\ \\ 
      z_1 - \tilde{\lambda} z_2 + z_4 &=& y_2 \, .
    \end{array}
\right.
\label{pfrv4}
\end{eqnarray}
\\
\nid
Assume a solution to (\ref{pfrv3}) has the form
\\
\begin{equation}
   z_n=c^{(1)}_n w^n_* + c^{(2)}_n (-w_*)^n +
  c^{(3)}_n w^{-n}_* + c^{(4)}_n (-w_*)^{-n} \, .
\label{pfrv5}
\end{equation}
\\
\nid
If
\\
\begin{eqnarray}
& &    \Delta c^{(1)}_n \ w^{n+1}_* + 
    \Delta c^{(2)}_n \ (-w_*)^{n+1} \nonumber \\  \nonumber \\
& & \ \ + \Delta c^{(3)}_n \ w^{-(n+1)}_* + 
    \Delta c^{(4)}_n \ (-w_*)^{-(n+1)} = 0 \, ,\label{pfrv6} \\ \nonumber \\
& &    \Delta c^{(1)}_n \ w^{n+2}_* + 
    \Delta c^{(2)}_n \ (-w_*)^{n+2}  \nonumber \\  \nonumber \\
& & \ \ + \Delta c^{(3)}_n \ w^{-(n+2)}_* + 
    \Delta c^{(4)}_n \ (-w_*)^{-(n+2)} = 0 \, , \label{pfrv7} \\ \nonumber \\
& &    \Delta c^{(1)}_n \ w^{n+3}_* + 
    \Delta c^{(2)}_n \ (-w_*)^{n+3} \nonumber \\  \nonumber \\
& & \ \ + \Delta c^{(3)}_n \ w^{-(n+3)}_* + 
    \Delta c^{(4)}_n \ (-w_*)^{-(n+3)} = 0 \, , \label{pfrv8} \\ \nonumber \\ 
& &  \Delta c^{(1)}_n \ w^{n+4}_* + 
    \Delta c^{(2)}_n \ (-w_*)^{n+4} \nonumber \\  \nonumber \\
& & \ \ + \Delta c^{(3)}_n \ w^{-(n+4)}_* + 
    \Delta c^{(4)}_n \ (-w_*)^{-(n+4)} = y_{n+2} \, ; \label{pfrv9}
\end{eqnarray}
\\
\nid
where $\Delta c^{(\ell)}_n = c^{(\ell)}_{n+1} - c^{(\ell)}_n$,
$(\ell = 1,2,3,4)$, then $z_n$ given in (\ref{pfrv5}) solves
(\ref{pfrv3}).  Solving (\ref{pfrv6}--\ref{pfrv9}), we have
\\
\begin{equation}
     \Delta c^{(\ell)}_n = (-1)^{\ell}\ y_{n+2} \ 
         \frac{D^{(\ell)}_n}{W_n} \, , 
         \quad (\ell = 1,2,3,4)
\label{pfrv10}
\end{equation}
\\
\nid
where
\\
\begin{eqnarray}
  W_n &=& \left|
    \begin{array}{cccc}
     w^{n+1}_* & (-w_*)^{n+1} & w^{-(n+1)}_* & (-w_*)^{-(n+1)} \\ \\ 
     w^{n+2}_* & (-w_*)^{n+2} & w^{-(n+2)}_* & (-w_*)^{-(n+2)} \\ \\ 
     w^{n+3}_* & (-w_*)^{n+3} & w^{-(n+3)}_* & (-w_*)^{-(n+3)} \\ \\ 
     w^{n+4}_* & (-w_*)^{n+4} & w^{-(n+4)}_* & (-w_*)^{-(n+4)} \\     
    \end{array}
         \right| \, , \label{pfrv11} \\ \nonumber \\ \nonumber \\
  D^{(1)}_n &=& \left|
    \begin{array}{ccc}
     (-w_*)^{n+1} & w^{-(n+1)}_* & (-w_*)^{-(n+1)} \\ \\ 
     (-w_*)^{n+2} & w^{-(n+2)}_* & (-w_*)^{-(n+2)} \\ \\ 
     (-w_*)^{n+3} & w^{-(n+3)}_* & (-w_*)^{-(n+3)} \\
    \end{array}
         \right|  \nonumber \\ \nonumber \\ \nonumber \\
&=& 2w^{-n}_* \left[ w^{-4}_* -1 \right] \, ,
\label{pfrv12} \\ \nonumber \\ \nonumber \\
  D^{(2)}_n &=& \left|
    \begin{array}{ccc}
     w_*^{n+1} & w^{-(n+1)}_* & (-w_*)^{-(n+1)} \\ \\
     w_*^{n+2} & w^{-(n+2)}_* & (-w_*)^{-(n+2)} \\ \\
     w_*^{n+3} & w^{-(n+3)}_* & (-w_*)^{-(n+3)} \\
    \end{array}
         \right|  \nonumber \\ \nonumber \\ \nonumber \\
&=& 2 (-w_*)^{-n} \left[ 1- w^{-4}_* \right] \, , \label{pfrv13}\\
\nonumber \\ \nonumber \\
  D^{(3)}_n &=& \left|
    \begin{array}{ccc}
     w_*^{n+1} & (-w_*)^{n+1} & (-w_*)^{-(n+1)} \\ \\
     w_*^{n+2} & (-w_*)^{n+2} & (-w_*)^{-(n+2)} \\ \\
     w_*^{n+3} & (-w_*)^{n+3} & (-w_*)^{-(n+3)} \\
    \end{array}
         \right|  \nonumber \\ \nonumber \\ \nonumber \\
&=& 2 w_*^n \left[ w^4_* -1 \right] \, , \label{pfrv14}\\ \nonumber \\
\nonumber \\
  D^{(4)}_n &=& \left|
    \begin{array}{ccc}
     w_*^{n+1} & (-w_*)^{n+1} & w_*^{-(n+1)} \\ \\ 
     w_*^{n+2} & (-w_*)^{n+2} & w_*^{-(n+2)} \\ \\
     w_*^{n+3} & (-w_*)^{n+3} & w_*^{-(n+3)} \\
    \end{array}
         \right|  \nonumber \\ \nonumber \\ \nonumber \\
&=& 2 (-w_*)^n \left[ 1- w^4_* \right] \, . \label{pfrv15}
\end{eqnarray}
\\
\nid
The $W_n$ defined in (\ref{pfrv11}) satisfies the Wronskian
relation
\\
\begin{equation}
  W_{n+1} = W_n \, .
\label{pfrv16}
\end{equation}
\\
\nid
The representation (\ref{pfrv10}) can be extended to $n \geq 0$.
Choose $c^{(\ell)}_0 = 0$, we have
\\
\begin{equation}
  c^{(\ell)}_n = \sum^{n-1}_{j=0} \Delta c^{(\ell)}_j = 
     \frac{1}{W_0} \ \sum^{n-1}_{j=0} D^{(\ell)}_j \ (-1)^{\ell}
     \ y_{j+2} \, , 
  \quad (n \geq 1) \, .
\label{pfrv17}
\end{equation}
\\
\nid
The expressions (\ref{pfrv12}--\ref{pfrv15}) lead to
\\
\begin{eqnarray}
  c^{(1)}_n &=& \frac{2 \left[ 1-w^{-4}_* \right]}{W_0}
              \ \sum^{n-1}_{j=0} w^{-j}_*\ y_{j+2} \, , \label{pfrv18} \\
\nonumber \\
  c^{(2)}_n &=& \frac{2 \left[ 1-w^{-4}_* \right]}{W_0}
              \ \sum^{n-1}_{j=0} (-w_*)^{-j} \ y_{j+2} \, ,  \label{pfrv19} \\
\nonumber \\
  c^{(3)}_n &=& \frac{2 \left[ 1-w^4_* \right]}{W_0}
              \ \sum^{n-1}_{j=0} w_*^{j} \ y_{j+2} \, , \label{pfrv20} \\
\nonumber \\
   c^{(4)}_n &=& \frac{2 \left[ 1-w^4_* \right]}{W_0}
              \ \sum^{n-1}_{j=0} (-w_*)^j \ y_{j+2} \, .\label{pfrv21}
\end{eqnarray}
\\
\nid
With these representations of $c^{(\ell)}_n$, $\ z_n$ given by
(\ref{pfrv5}) is a special solution to (\ref{pfrv3}).  The
general solution to (\ref{pfrv3}) is
\\
\begin{eqnarray}
z_n &=& 
      \left( c^{(1)}_n + a^{(1)} \right) w^n_* +
      \left( c^{(2)}_n + a^{(2)} \right) (-w_*)^n 
      \nonumber   \\ \nonumber   \\
     && \ + \left( c^{(3)}_n + a^{(3)} \right) w^{-n}_* +
        \left( c^{(4)}_n + a^{(4)} \right) (-w_*)^{-n} \, , 
  \label{pfrv22}
\end{eqnarray}
\\
\nid
where $a^{(\ell)}$ $(\ell = 1,2,3,4)$ are arbitrary constants.
Set
\\
\begin{eqnarray}
  a^{(3)} &=& - \ \frac{2 \left[ 1-w^4_* \right]}{W_0}  
            \ \sum^{\infty}_{j=0} w^j_* \ y_{j+2} \, , \label{pfrv23} \\
\nonumber \\
   a^{(4)} &=& - \ \frac{2 \left[ 1-w^4_* \right]}{W_0} \  
            \sum^{\infty}_{j=0} (-w_*)^j \ y_{j+2} \, . \label{pfrv24}
\end{eqnarray}
\\
\nid
Let
\\
\begin{equation}
   f_n = c^{(1)}_n w_*^n +  c^{(2)}_n (-w_*)^n + 
       (c^{(3)}_n + a^{(3)}) w_*^{-n} + (c^{(4)}_n + a^{(4)})
       (-w_*)^{-n} \, ;
\label{pfrv25}
\end{equation}
\\
\nid
then from expressions (\ref{pfrv18}--\ref{pfrv21}) and 
(\ref{pfrv23}--\ref{pfrv25}), we have
\\
\begin{eqnarray}
 f_n &=& \frac{2 \left[ 1-w_*^{-4} \right]}{W_0} 
       \ \sum^{n-1}_{j=0} \left[ w_*^{n-j} + (-w_*)^{n-j} \right]
        y_{j+2} \nonumber \\ \nonumber \\
     & & - \frac{2 \left[ 1-w_*^4 \right]}{W_0} 
       \sum^{\infty}_{j=n} \left[ w_*^{j-n} + (-w_*)^{j-n} \right]
        y_{j+2} \, .  \label{pfrv26}
\end{eqnarray}
\\
\nid
Finally,
\\
\begin{equation}
 z_n = a^{(1)} w^n_* + a^{(2)} (-w_*)^n + f_n \, . 
\label{pfrv27}
\end{equation}
\\
\nid
Next we choose $a^{(1)}$ and $a^{(2)}$ to satisfy the constraints
(\ref{pfrv4}):
\\
\begin{equation}
 M \left(
   \begin{array}{c}
         a^{(1)} \\
         a^{(2)} 
   \end{array}
\right) = 
    \left[ 
   \begin{array}{c}
         y_1 + \tilde{\lambda} f_1 - f_2 -f_3 \\[1ex]
         y_2 - f_1 + \tilde{\lambda} f_2 - f_4 
   \end{array}
\right] \, ,
\label{pfrv28}
\end{equation}
\\
\nid
where
\\
\begin{eqnarray*}
 M= \left[
   \begin{array}{cc}
       - \tilde{\lambda} w_* + w^2_* + w^3_*   &
        \tilde{\lambda} w_* + w_*^2 - w_*^3 \\ \\
        w_* - \tilde{\lambda} w^2_* + w^4_* &
        - w_* - \tilde{\lambda} w^2_* + w^4_*
   \end{array}
    \right] \, .  
\end{eqnarray*}
\\
\nid
As shown in the proof of Theorem~\ref{pspec}, $M$ is
nonsingular.  Then,
\\
\begin{eqnarray}
 \left(
   \begin{array}{c}
         a^{(1)} \\
         a^{(2)} 
   \end{array}
\right) = M^{-1}
    \left[ 
   \begin{array}{c}
         y_1 + \tilde{\lambda} f_1 - f_2 -f_3 \\[1ex]
         y_2 - f_1 + \tilde{\lambda} f_2 - f_4 
   \end{array}
\right] \, , 
\label{pfrv29}
\end{eqnarray}
\\
\nid
where
\\
\begin{displaymath}
  M^{-1} = \tilde{\k}^{-1}
           \left[
              \begin{array}{cc}
                -w_* - \tilde{\lambda} w^2_* + w^4_* &
                - \tilde{\lambda} w_* - w^2_* + w^3_* \\ \\
                -w_* + \tilde{\lambda} w^2_* - w^4_* &
                - \tilde{\lambda} w_* + w^2_* + w^3_* \\
              \end{array}
            \right] \, ,
\end{displaymath}
\\
\nid
where $\tilde{\k} = 2w_*^3 \left( w^4_* - 2 \tilde{\lambda} w^2_* +
\tilde{\lambda}^2 -1 \right)$.
Thus
\\
\begin{equation}
  z_n = \bigg ( w^n_* \, , (-w_*)^n \bigg ) \ M^{-1}\ 
        \left[
          \begin{array}{c}
            y_1 + \tilde{\lambda} f_1 - f_2 - f_3 \\[1ex]
            y_2 - f_1 + \tilde{\lambda} f_2 - f_4
          \end{array}
        \right] + f_n
\label{pfrv30}
\end{equation}
\\
\nid
solves Equation (\ref{pfrv2}).  Rewrite $f_n$ given in
(\ref{pfrv26}) as follows:
\\
\begin{equation}
    f_n = \sum^{\infty}_{j=1} g(n,j)\  y_j \, , \label{pfrv31}
\end{equation}
\\
\nid
where
\\
\begin{eqnarray*}
  g(n,j) = \left\{
      \begin{array}{cl}
          0 \, , \quad & j=1 \, ; \\ \\
          \frac{2[1-w_*^{-4}]}{W_0} \ \left[ w^{n-j+2}_* +
                (-w_*)^{n-j+2} \right] \, , \quad &
                2 \leq j \leq n+1 \, ; \\ \\ [1ex]
          - \ \frac{2[1-w^4_*]}{W_0}\  \left[ w^{j-n-2}_* + 
                (-w_*)^{j-n-2} \right] \, , \quad &
                j \geq n+2 \, .
      \end{array}
      \right.
\end{eqnarray*}
\\
\nid
Rewrite $z_n$ given in (\ref{pfrv30}) as follows:
\\
\begin{equation}
    z_n= \sum^{\infty}_{j=1} G(n,j) \ y_j \, , \label{pfrv32}
\end{equation}
\\
\nid
where
\\
\begin{eqnarray}
  G(n,j) &=& \bigg ( w^n_* , \ (-w_*)^n \bigg )\ M^{-1} \nonumber \\ 
\nonumber \\       
& &\left(
        \begin{array}{c}
          \delta_{1,j} + \tilde{\lambda}\ g(1,j) -
            g(2,j) - g(3,j) \\[1ex]
          \delta_{2,j} -  g(1,j) + \tilde{\lambda}\ 
            g(2,j) - g(4,j) 
        \end{array}
       \right) + g(n,j) \, ,
\label{pfrv33}
\end{eqnarray}
\\
\nid
where $\delta_{\ell,j}$ is the Kronecker delta:
$\delta_{\ell,j} =1$ $(\ell=j)$, $\delta_{\ell,j} =0$ $(\ell
\neq j)$.
From the expression (\ref{pfrv33}), we see that there exists a
constant $K$ independent of $n, j$; such that
\begin{eqnarray}
  \sum^{\infty}_{j=1} \left| G(n,j) \right| &\leq& K \, , \quad
      \forall \ n=1,2, \cdots \, ;  \label{pfrv34} \\ \nonumber \\
  \sum^{\infty}_{n=1} \left| G(n,j) \right| &\leq& K \, , \quad
      \forall \ j=1,2, \cdots \, .  \label{pfrv35}
\end{eqnarray}
\\
\nid
Then, we have the $\ell_{\infty}$ norm relation,
\\
\begin{eqnarray}
  \| z\|_{\infty} = \sup_n \left| z_n \right| 
         & \leq & \sup_n \ \sum^{\infty}_{j=1} 
                    \left| G(n,j)  \right| \ \left| y_j \right| 
                    \nonumber\\  \nonumber\\
         & \leq & \left[ \sup_n \ \sum^{\infty}_{j=1} 
                    \left| G(n,j) \right| \right] \ \| y \|_{\infty} 
                  \nonumber\\  \nonumber\\
         & \leq & K \ \| y \|_{\infty} \, ,\label{pfrv36}
\end{eqnarray}
\\
\nid
and the $\ell_1$ norm relation,
\\
\begin{eqnarray}
  \| z \|_1 = \lim_{N \to \infty} \  
              \left[ \sum^N_{n=1}  \ \left| z_n \right| \right]
       & \leq & \lim_{N \to \infty} \left[ \sum^N_{n=1} \sum^{\infty}_{j=1} 
                    \left| G(n,j)\right| \ \left| y_j \right|
                  \right] \nonumber \\ \nonumber \\
       & = & \lim_{N \to \infty} \left[ \sum^{\infty}_{j=1} 
                    \left|  y_j \right| \sum^N_{n=1} \left|
                      G(n,j) \right| \right] \nonumber \\ \nonumber \\
       & \leq & K \ \sum^{\infty}_{j=1} \left| y_j \right| = K \ \| y \|_1 \, .
\label{pfrv37}
\end{eqnarray}
\\
Thus the linear operator defined in (\ref{pfrv32}) which maps $y$ into $z$,
is bounded in $\ell_{\infty}$ and $\ell_1$.  Therefore, by Riesz
convexity theorem \cite{Rie26} \cite{CZ50} \cite{Dur60}, $(\LL_{\tilde{B}} -
\tilde{\lambda}I)^{-1}$ defined in (\ref{pfrv32}) is bounded in
$\ell_2$.  Since by Theorems~\ref{pspec} and \ref{rspec},
$(\LL_{\tilde{B}} - \tilde{\lambda}I)^{-1}$ exists and is
everywhere densely defined and is bounded, we have
$\tilde{\lambda} \in \rho(\tilde{B})$.  In summary, we have
shown that if $\tilde{\lambda} \in C_{\tilde{B}}$, then
$\tilde{\lambda} \in \sigma_c(\tilde{B})$; and if
  $\tilde{\lambda} \notin C_{\tilde{B}}$, then $\tilde{\lambda}
    \in \rho(\tilde{B})$; thus, $\sigma_c(\tilde{B}) =
    C_{\tilde{B}}$ and $\rho(\tilde{B}) = (C_{\tilde{B}})'$.
    Equivalently, $\sigma_c(B) = C_B$ and $\rho(B) = (C_B)'$.
$\Box$

In summary, the spectrum of $\LL_B$ is as depicted in Figure \ref{splb}.

\subsection{The Spectra of the Linear Operator $\LL_A$}

\begin{figure}[ht]
  \begin{center}
    \leavevmode
      \setlength{\unitlength}{2ex}
  \begin{picture}(36,27.8)(-18,-12)
    \thicklines
\put(0,-14){\vector(0,1){28}}
\put(-18,0){\vector(1,0){36}}
\put(0,15){\makebox(0,0){$\Im \{ \la \}$}}
\put(18.5,0){\makebox(0,0)[l]{$\Re \{ \la \}$}}
\put(0.1,-7){\line(0,1){14}}
\put(.2,-.2){\makebox(0,0)[tl]{$0$}}
\put(-0.2,-7){\line(1,0){0.4}}
\put(-0.2,7){\line(1,0){0.4}}
\put(2.0,-6.4){\makebox(0,0)[t]{$-i2|b|$}}
\put(2.0,7.6){\makebox(0,0)[t]{$i2|b|$}}
\end{picture}
  \end{center}
\caption{The continuous spectrum of $\LL_B$ and $\LL_A$.}
\label{splb}
\end{figure}
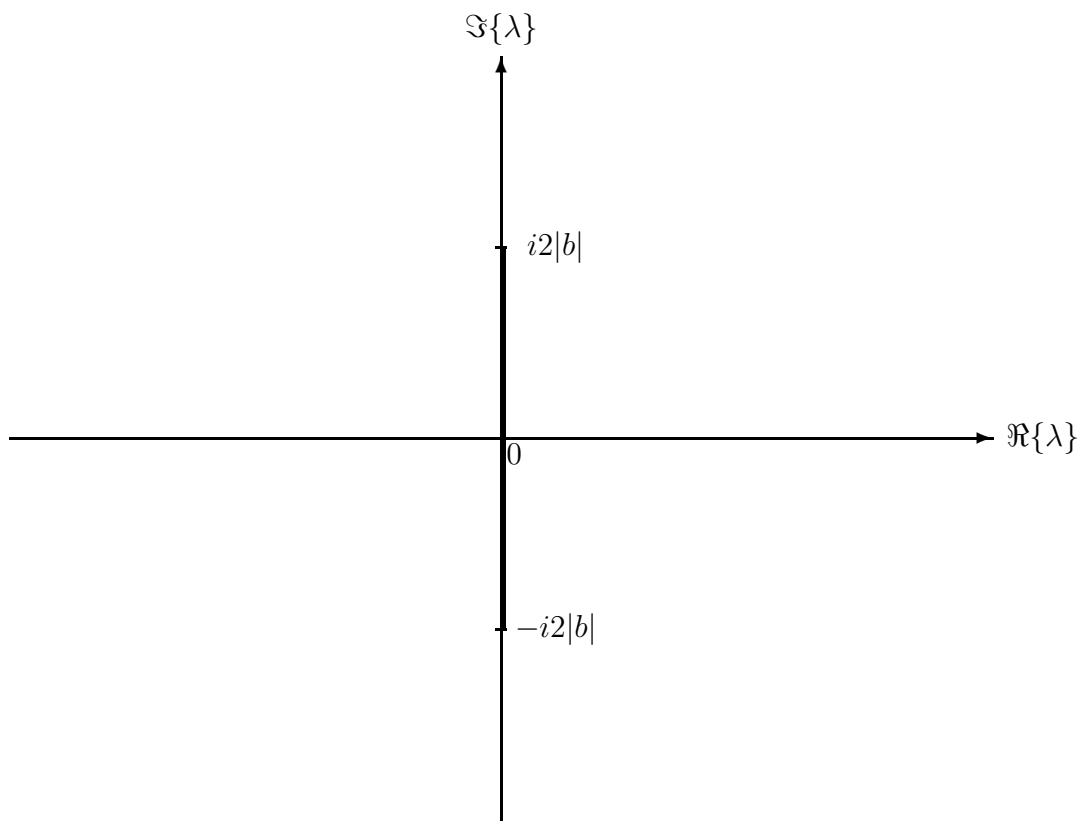

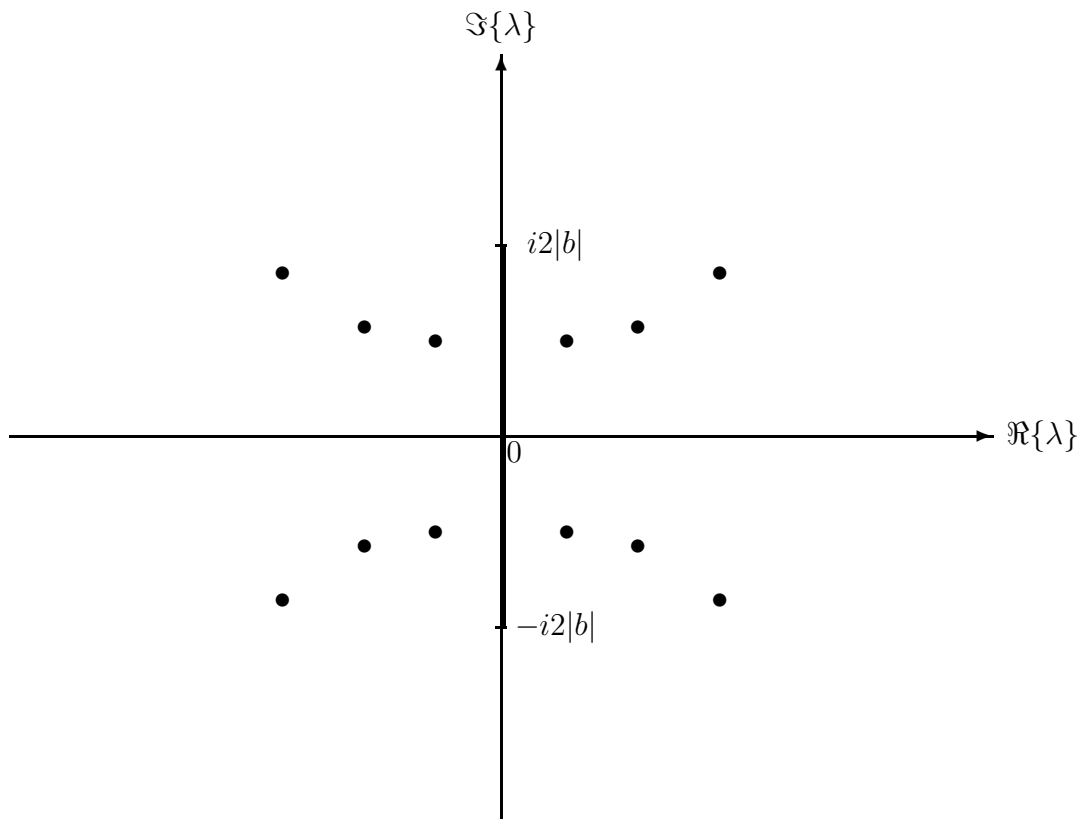
\begin{figure}[ht]
  \begin{center}
    \leavevmode
      \setlength{\unitlength}{2ex}
  \begin{picture}(36,27.8)(-18,-12)
    \thicklines
\put(0,-14){\vector(0,1){28}}
\put(-18,0){\vector(1,0){36}}
\put(0,15){\makebox(0,0){$\Im \{ \la \}$}}
\put(18.5,0){\makebox(0,0)[l]{$\Re \{ \la \}$}}
\put(0.1,-7){\line(0,1){14}}
\put(.2,-.2){\makebox(0,0)[tl]{$0$}}
\put(-0.2,-7){\line(1,0){0.4}}
\put(-0.2,7){\line(1,0){0.4}}
\put(2.0,-6.4){\makebox(0,0)[t]{$-i2|b|$}}
\put(2.0,7.6){\makebox(0,0)[t]{$i2|b|$}}
\put(2.4,3.5){\circle*{0.5}}
\put(-2.4,3.5){\circle*{0.5}}
\put(2.4,-3.5){\circle*{0.5}}
\put(-2.4,-3.5){\circle*{0.5}}
\put(5,4){\circle*{0.5}}
\put(-5,4){\circle*{0.5}}
\put(5,-4){\circle*{0.5}}
\put(-5,-4){\circle*{0.5}}
\put(8,6){\circle*{0.5}}
\put(-8,6){\circle*{0.5}}
\put(8,-6){\circle*{0.5}}
\put(-8,-6){\circle*{0.5}}
\end{picture}
  \end{center}
\caption{The spectrum of $\LL_A$ in case 2.}
\label{spla2}
\end{figure}

Now we apply Weyl's essential spectrum theorem \cite{RS78} to obtain 
the spectral theorem for $\LL_A$.
\begin{theorem}[The Spectral Theorem of $\LL_A$]
\begin{enumerate}
\item If $\Sg_{\hat{k}} \cap \bar{D}_{|p|} = \emptyset$, then the entire
$\ell_2$ spectrum of the linear operator $\LL_A$ is its continuous spectrum 
which is the spectral curve $C_B$ defined in (\ref{srb}), i.e. $\sg (\LL_A) 
= \sg_c (\LL_A) = C_B$. See Figure \ref{splb}.
\item If $\Sg_{\hat{k}} \cap \bar{D}_{|p|} \neq \emptyset$, then the entire
essential $\ell_2$ spectrum of the linear operator $\LL_A$ is its 
continuous spectrum 
which is the spectral curve $C_B$ defined in (\ref{srb}), i.e. 
$\sg_{\mbox{ess}} (\LL_A) = \sg_c (\LL_A) = C_B$. That is, the residual 
spectrum of $\LL_A$ is empty, $\sg_r (\LL_A) = \emptyset$. The point spectrum 
of $\LL_A$ is symmetric with respect to both real and imaginary axes. 
See Figure \ref{spla2}.
\end{enumerate}
\label{spthla}
\end{theorem}

Proof: First, we want to show that in both cases, the residual spectrum of 
$\LL_A$ is empty. By Weyl's essential spectrum theorem \cite{RS78}, the 
essential spectrum of $i\LL_A$ is the same with the essential spectrum of 
$i\LL_B$, $\sg_{\mbox{ess}} (i\LL_A) = \sg_{\mbox{ess}} (i\LL_B) =iC_B$. 
Let $i \la_r \in \sg_r (i\LL_A)$, then $i \la_r \in iC_B$. By the argument 
in the proof of Theorem \ref{rspec}, $i \la_r \in \sg_p ((i\LL_A)^*)$, 
where $(i\LL_A)^*$ is the adjoint of $i\LL_A$,
\[
(i\LL_A)^* = i\LL_B + (i\LL_C)^*\ .
\]
By Weyl's essential spectrum theorem \cite{RS78}, the essential spectrum 
of $(i\LL_A)^*$ is the same with the essential spectrum of 
$i\LL_B$, $\sg_{\mbox{ess}} ((i\LL_A)^*) = \sg_{\mbox{ess}} (i\LL_B) =iC_B$.
Thus, $i \la_r \in \sg_{\mbox{ess}} ((i\LL_A)^*)$. Since $\sg_{\mbox{ess}} 
((i\LL_A)^*)$ and $\sg_p ((i\LL_A)^*)$ are disjoint, $\sg_r (i\LL_A) = \emptyset$. The claim $\sg_p (\LL_A) = \emptyset$ in case 1 follows fom the proof of 
Lemma \ref{lenev} and the fact that the spectral curve $C_B$ corresponds 
to $\mbox{Re}\{ \ta \}=0$ and $|\ta| \leq 2$. The property of $\sg_p (\LL_A)$
in case 2 has been proved in Theorem \ref{prpev}. Then Weyl's essential 
spectrum theorem implies the rest of the claims. $\Box$
\begin{remark}
By the above theorem, the computation of eigenvalues is reduced to the 
case that $\Sg_{\hat{k}} \cap \bar{D}_{|p|} \neq \emptyset$. By Corollary  
\ref{cgtle}, if $\la \not \in C_B = \sg_{\mbox{ess}} (\LL_A)$, the two 
continued fractions (\ref{cfr4}) and (\ref{cfr6}) converge. And solutions 
of equation (\ref{cfr17}) lead to eigenvalues.
\end{remark}
\begin{remark}
The width of the continuous spectrum $\sg_c (\LL_A)$ is $4|b|$, where 
$b=-a|p|^{-2}$ and $a = \frac{1}{2} \left| \Gamma \right| 
\left| \begin{array}{cc}
p_1 & \hat{k}_1 \\
p_2 & \hat{k}_2
  \end{array}
\right| \, 
$. Although $|a|$ can increase to infinity as $|k|$ increases to infinity,
$a$ is essentially a scaling factor for $\LL_A$ as can be seen in the 
expression for the infinite-matrix $A$.
\end{remark}
Next we discuss an alternative way of representing eigenvalues. This approach 
is not useful for practical computation.
Consider the linear difference equation,
\\
\begin{equation}
(A- \lambda I)\ z = 0 \, ,
\label{las1}
\end{equation}
\\
\nid
where $A$ defined in (\ref{mata}) is the representation matrix of
$\LL_A$.  Explicitly,
\\
\begin{eqnarray}
  \left\{
      \begin{array}{ll}
        \rho_{n-1} \ z_{2(n-1)} - \hat{\lambda} \ z_{2n} +
        \rho_{n+1} \ z_{2(n+1)} = 0 \, , \quad & (n \geq 2), \\ \\
        \rho_{-n+1} \ z_{2(n-1)+1} - \hat{\lambda} \ z_{2n+1} +
        \rho_{-n-1} \ z_{2(n+1)+1} = 0 \, , \quad & (n \geq 1),
      \end{array}
  \right.
\label{las2}
\end{eqnarray}
\\
\nid
under the constraints
\\
\begin{eqnarray}
   \left\{
    \begin{array}{l}
      - \hat{\lambda} \ z_1 + \rho_1 \ z_2 + \rho_{-1}\ z_3 = 0 \, , \\ \\
        \rho_0 \ z_1 - \hat{\lambda}\ z_2 + \rho_2 \ z_4 = 0 \, ,
    \end{array}
  \right.
\label{las3}
\end{eqnarray}
\\
\nid
where $\hat{\lambda} = (ia)^{-1}\ \lambda$.  By Theorem
\ref{cpt}, the linear operator $\LL_A$ is a compact perturbation
of $\LL_B$.  Thus, the difference equation (\ref{las2}) is of
Poincar\'{e}-Perron type.
The Poincar\'{e}-Perron theorem stated specifically for the
difference equation (\ref{las1}) is as follows \cite{Per10a},
\cite{Per10b}, \cite{Dur60}:
\begin{theorem}[Poincar\'{e}-Perron Theorem]
For any $\lambda \in C$, let $w_*$, $-\ w_*$, $\frac{1}{w_*}$, $-\ 
\frac{1}{w_*}$ be the roots of the characteristic polynomial $f_B
(w, \lambda)$ defined in (\ref{cpl1}), which are given in
(\ref{cpl4}), where $\left| w_* \right| \leq 1$.  Then there
exists a fundamental set of solutions $z^{(j)}_n$ $(j=1,2,3,4)$ to
the difference equation (\ref{las2}), such that
\\
\begin{eqnarray}
  \limsup_{n \to \infty}\ 
      \left| z^{(j)}_n \right|^{\frac{1}{n}} &=& \left| w_* \right|, 
      \quad (j=1,2) \, , \nonumber \\ \label{las4} \\
\limsup_{n \to \infty} \ \left| z^{(j)}_n
      \right|^{\frac{1}{n}} &=& \frac{1}{\left| w_* \right|}, 
      \quad (j=3,4) \, . \nonumber
\end{eqnarray}
\label{PPT}
\end{theorem}
It is easy to see that 
\\
\[
z^{(j)} \in \ell_2,\  (j=1,2); \ \ z^{(j)}
\notin \ell_2,\  (j=3,4); 
\]
\\
\nid
if $\left| w_* \right| < 1$.  By
definition, when $\lambda \notin C_B$ (defined in (\ref{srb})),
$\left| w_* \right| <1$.  Next we study the conditions for the
point spectrum of $\LL_A$.  Let
\\
\begin{equation}
  z_n = c_1 \ z^{(1)}_n + c_2 \ z^{(2)}_n \, .
\label{las5}
\end{equation}
\\
\nid
Substitute $z_n$ into the constraints (\ref{las3}), we have
\\
\begin{equation}
   M \left(
    \begin{array}{c}
      c_1 \\
      c_2
    \end{array}
    \right) = 0 \, ,
\label{las6}
\end{equation}
\\
\nid
where
\\
\begin{eqnarray}
     M= \left(
       \begin{array}{cc}
         -\hat{\lambda} z_1^{(1)} + \rho_1 z_2^{(1)} +
              \rho_{-1} z_3^{(1)} &
          -\hat{\lambda} z_1^{(2)} + \rho_1 z_2^{(2)} +
              \rho_{-1} z_3^{(2)} \\ \\
          \rho_0 z_1^{(1)} - \hat{\lambda} z_2^{(1)} +
              \rho_2 z_4^{(1)}  &
           \rho_0 z_1^{(2)} - \hat{\lambda} z_2^{(2)} +
              \rho_2 z_4^{(2)} 
       \end{array}
         \right) \, .
\label{las7}
\end{eqnarray}
\\
\begin{theorem}
If $\lambda \notin C_B$ (the spectral curve for $\LL_B$, defined
in (\ref{srb})), and $\det M=0$ (where $M$ is defined in
(\ref{las7})), then $\lambda \in \sigma_p(A)$ (the point
spectrum of $\LL_A$).
\label{ptspa1}
\end{theorem}

Proof:
  If $\det M=0$, then there is a nontrivial solution to
  (\ref{las6}).  Thus there is a nonzero solution to
  (\ref{las2}), which satisfies the constraints (\ref{las3}).
  Therefore, $\lambda$ is an eigenvalue. $\Box$

\newpage
\eqnsection{Conclusion}

In this paper, we study the linearized two-dimensional Euler equation 
at a stationary state. 
This equation decouples into infinite many invariant subsystems. Each 
invariant subsystem 
is shown to be a linear Hamiltonian system of infinite dimensions. 
Another important invariant besides the Hamiltonian for each invariant
subsystem is found, and is utilized to prove an ``unstable disk 
theorem'' through a simple Energy-Casimir argument. 
The eigenvalues of the 
linear Hamiltonian system are of four types: real pairs ($c,-c$), purely 
imaginary pairs ($id,-id$), quadruples ($\pm c\pm id$), and zero eigenvalues. 
The eigenvalues are studied through continued fractions.
The spectral equation for each invariant subsystem
is a Poincar\'{e}-type difference equation, i.e. it can be 
represented as the spectral equation of an infinite matrix operator, 
and the infinite matrix operator is a sum of a constant-coefficient 
infinite matrix operator and a compact infinite matrix operator. 
We have a complete spectral theory. The essential spectrum of each invariant 
subsystem is a bounded band of continuous spectrum. The point spectrum can 
be computed through continued fractions.

This study is the first step toward understanding the unstable manifold 
structures of stationary states of the two-dimensional Euler equation, 
which we believe to be the key for understanding two-dimensional turbulence.
In particular, we will be interested in investigating whether or not 
the unstable manifolds of 2D Euler equations are degenerate (i.e. figure 
eight structures). Degeneracy will imply that the dynamics of 2D Euler 
equations is not turbulent.

{\bf Acknowledgment:} This work was started at MIT, continued at 
Institute for Advanced Study, and finally completed at University 
of Missouri. The author had benefited a lot from discussions with 
Professor Thomas Witelski at MIT.

\newpage
\bibliography{leu}

\newpage

\newpage

\newpage

\newpage

\newpage

\newpage 

\newpage

Figure 1. Caption: An illustration on the locations of the modes ($k'=rk$) and 
($|k'|=|k|$) in the definitions of $E^1_k$ and $E^2_k$ 
(Proposition \ref{eman}).

Figure 2. Caption: An illustration of the classes $\Sg_{\hk}$ and the disk 
$\bar{D}_{|p|}$.

Figure 3. Caption: The quadruple of eigenvalues determined by equation (\ref{cfr17})
for the system led by the class $\Sg_{\hat{k}}$ labeled by $\hk = (1,0)^T$, 
when $p=(1,1)^T$.

Figure 4. Caption: The spectral curve $C_{\tilde{B}}\ $.

Figure 5. Caption: The continuous spectrum of $\LL_B$ and $\LL_A$.

Figure 6. Caption: The spectrum of $\LL_A$ in case 2.

\end{document}